%% file: main.tex
\documentclass[11pt,reqno]{amsart}
\input{96_standpacks}
\input{97_configs}
\input{98_tikz}

\usetikzlibrary{matrix,fit}
\usepackage{amssymb}
\usepackage{pdflscape}
\usepackage{amsmath}
\usepackage{forest}
\usepackage{float}
\usetikzlibrary{shapes.geometric}
\usetikzlibrary{fit,backgrounds,decorations.pathmorphing}
\usetikzlibrary{arrows.meta,decorations.pathreplacing}
\newcommand{\circled}[1]{\tikz[baseline=(char.base)]{\node[draw, circle, inner sep=0.3pt, line width=0.4pt] (char) {\scriptsize $#1$};}}
\newcommand{\Iso}{{\mathbin{\circled{\mathcal I}}}}
\newcommand{\Uni}{{\mathbin{\circled{\mathcal U}}}}

\forestset{
  circ/.style={
    draw,
    circle,
    inner sep=1pt,
    minimum size=4mm,
    font=\scriptsize
  }
}

\title[Benign Nonconvexity of Synchronization Landscape Induced by Graph Skeletons]{Benign Nonconvexity of Synchronization Landscape Induced by Graph Skeletons}
\author[Wu]{Hongjin Wu\textsuperscript{\normalfont *}}
\thanks{\textsuperscript{\normalfont *}%
ETH Z\"urich, Switzerland.
\texttt{hongjin-wu@outlook.com}.
This work was carried out while the first author was
a doctoral researcher at ETH Z\"urich. Current address: Centrum Wiskunde \& Informatica (CWI), The Netherlands.}

\author[Brandes]{Ulrik Brandes\textsuperscript{\normalfont\textdagger}}
\thanks{\textsuperscript{\normalfont\textdagger}%
ETH Z\"urich, Switzerland.
\texttt{ubrandes@ethz.ch}.}
\begin{document}
\begin{abstract}
We study the homogeneous Kuramoto model on a graph and
the associated nonconvex optimization problem
$\min_{\boldsymbol{\theta}\in\mathbb{R}^n}
-\frac12\sum_{1\leq i,j\leq n}A_{ij}\cos(\theta_i-\theta_j)$.
The objective defines an energy over configurations of
points on the unit circle and serves as a Lyapunov potential
for the dynamics.
We prove that every connected quasi-threshold graph is
second-order globally synchronizing: every second-order
stationary point is a global minimizer corresponding to
full synchronization.
Consequently, the dynamics converge to full synchronization
from almost every initial condition.
These graphs are precisely the comparability graphs of
partially ordered sets induced by rooted trees.
Viewing these trees as graph skeletons, we establish an
upward propagation mechanism: synchronization within each
child subtree forces the parent and all its descendants
to synchronize.
The argument proceeds from the leaves to the root and
relies on phasor geometry and second-order optimality
conditions.
Our result provides a structural route to global
synchronization that complements existing results based
on minimum-degree or spectral proximity to complete graphs.
\end{abstract}

\maketitle
\tableofcontents
\input{01_introduction}
\input{02_quasi}
\subsection*{Acknowledgements}
We would like to thank Afonso S. Bandeira for insightful discussions, and Monique Laurent for valuable feedback on our work that helped improve this manuscript.
\bibliographystyle{alpha}
\bibliography{ref}
\addtocontents{toc}{\protect\setcounter{tocdepth}{1}}%
\setcounter{tocdepth}{1}%
\vspace*{1cm}

\end{document}

%% file: 96_standpacks.tex
\usepackage[utf8]{inputenc}
\usepackage[english]{babel}
\usepackage{ifthen} 

\usepackage{amsmath} 
\usepackage{amssymb} 
\usepackage{mathrsfs} 
\usepackage{mathtools} 
\usepackage{bbm} 
\usepackage{cancel} 
\usepackage{mathdots} 
\usepackage{xfrac} 

\usepackage{amsthm} 
\usepackage[shortlabels]{enumitem}

\usepackage{lipsum} 
\usepackage{verbatim} 
\usepackage{pifont} 
\usepackage[dvipsnames]{xcolor} 
\usepackage{lmodern} 
\usepackage{marginnote} 

\usepackage{wrapfig} 
\usepackage{caption} 

\usepackage{array} 
\usepackage{multirow} 
\usepackage{multicol} 
\usepackage{colortbl} 
\usepackage{makecell} 
\usepackage{arydshln} 

\usepackage{tikz} 
\usetikzlibrary{automata,positioning} 
\usepackage{tkz-euclide} 
\usepackage[customcolors,shade]{hf-tikz} 

\usepackage{geometry} 
\usepackage{titleps} 
\usepackage{footnote} 

\usepackage{hyperref} 
\usepackage{cleveref} 


%% file: 97_configs.tex
\newtheorem{theorem}{Theorem}[section]

\newtheorem{lemma}[theorem]{Lemma}
\newtheorem{corollary}[theorem]{Corollary}

\theoremstyle{definition}
\newtheorem{definition}[theorem]{Definition}

\newcolumntype{L}[1]{>{\raggedright\let\newline\\\arraybackslash\hspace{0pt}}m{#1}}
\newcolumntype{C}[1]{>{\centering\let\newline\\\arraybackslash\hspace{0pt}}m{#1}}
\newcolumntype{R}[1]{>{\raggedleft\let\newline\\\arraybackslash\hspace{0pt}}m{#1}}

\usepackage{float}

\hypersetup{
    colorlinks=true,
    linkcolor=blue,
    filecolor=magenta,
    urlcolor=red,
    citecolor=blue,
}


%% file: 98_tikz.tex
\usetikzlibrary{decorations.pathmorphing}
\usetikzlibrary{decorations.markings}
\usetikzlibrary{decorations.text}

\usetikzlibrary{fit}
\usetikzlibrary{backgrounds}
\usepackage{subcaption}
\usepackage{adjustbox} 

\definecolor{cfefefe}{RGB}{254,254,254}
\definecolor{c0b0b0b}{RGB}{11,11,11}
\definecolor{cc7c6c3}{RGB}{199,198,195}
\definecolor{cdfdedf}{RGB}{223,222,223}
\definecolor{c0a0a09}{RGB}{10,10,9}
\definecolor{c0a0a0a}{RGB}{10,10,10}
\definecolor{c080808}{RGB}{8,8,8}
\definecolor{cbfbcb6}{RGB}{191,188,182}
\definecolor{cc1beb7}{RGB}{193,190,183}
\definecolor{cb5b4b1}{RGB}{181,180,177}
\definecolor{cd8d8d8}{RGB}{216,216,216}
\definecolor{cb8b7b3}{RGB}{184,183,179}
\definecolor{cc4c2c0}{RGB}{196,194,192}
\definecolor{c1b1b1a}{RGB}{27,27,26}
\definecolor{c9e9991}{RGB}{158,153,145}
\definecolor{cbdbab4}{RGB}{189,186,180}
\definecolor{cb9b7b4}{RGB}{185,183,180}
\definecolor{cb9b8b5}{RGB}{185,184,181}
\definecolor{c1a1a19}{RGB}{26,26,25}
\definecolor{ccfcecb}{RGB}{207,206,203}
\definecolor{cb3aea7}{RGB}{179,174,167}
\definecolor{cc9c7c1}{RGB}{201,199,193}
\definecolor{c6b6a68}{RGB}{107,106,104}
\definecolor{cd3d2cf}{RGB}{211,210,207}
\definecolor{c8b8b8b}{RGB}{139,139,139}
\definecolor{ca6a5a2}{RGB}{166,165,162}
\definecolor{caaa9a6}{RGB}{170,169,166}
\definecolor{c6e6c67}{RGB}{110,108,103}
\definecolor{ca5a4a1}{RGB}{165,164,161}
\definecolor{c9d9c9a}{RGB}{157,156,154}
\definecolor{cb8b7b4}{RGB}{184,183,180}
\definecolor{c1c1c1b}{RGB}{28,28,27}
\definecolor{ca39e96}{RGB}{163,158,150}
\definecolor{ca8a6a0}{RGB}{168,166,160}
\definecolor{cbcbbb8}{RGB}{188,187,184}
\definecolor{c505050}{RGB}{80,80,80}
\definecolor{c9f9b93}{RGB}{159,155,147}
\definecolor{cb4b3b0}{RGB}{180,179,176}
\definecolor{c2c2c2b}{RGB}{44,44,43}
\definecolor{ce2e2e3}{RGB}{226,226,227}
\definecolor{ccbc9c3}{RGB}{203,201,195}

%% file: 01_introduction.tex
\section{Introduction}\label{sec:intro}
Nonconvex optimization is generally computationally hard, yet some problems exhibit surprisingly simple landscapes.
A prominent example is the Burer--Monteiro factorization of semidefinite programming \cite{burer2003nonlinear,boumal2020deterministic}, where, under suitable conditions, even simple local algorithms exhibit strong practical performance. A theoretical explanation for this success is that, in such cases, every second-order stationary point is globally optimal, a phenomenon often referred to as benign nonconvexity \cite{GeLeeMa2017,SunQuWright2018,BoumalAbsilCartis2018}. Understanding when this property holds is therefore crucial for distinguishing problems where local search can reliably reach a global optimum from those where it may fail.

The synchronization problem on graphs provides a particularly clean setting for studying when such benign landscapes arise.
Let $G=(V,E)$ be a connected graph with adjacency matrix $A$. Each vertex $i\in V$ is assigned a unit vector $\boldsymbol v_i\in\mathbb{S}^1$.
The goal is to find configurations $(\boldsymbol v_1, \ldots, \boldsymbol v_n) \in (\mathbb S^1)^n$ that minimize the following energy function:
\begin{equation}\label{E_1}
\min_{\forall i \in [n], \ v_i\in \mathbb{R}^2,\ \|\boldsymbol v_i\|=1}
E_G(\boldsymbol v_1,\ldots,\boldsymbol v_n)
\;=\;
\frac{1}{2}\sum_{1 \le i,j \le n} A_{ij} \bigl(1 - \langle\boldsymbol  v_i, \boldsymbol v_j \rangle\bigr).
\end{equation}
Here, each pair of adjacent vertices $(i,j)$ prefers to synchronize, meaning that configurations with $\boldsymbol v_i = \boldsymbol v_j$ are energetically favorable, while discrepancies are energetically unfavorable.
Despite its simple form, this energy is nonconvex in general and admits spurious second-order stationary points.
However, under certain graph structures, the landscape becomes remarkably simple: every second-order stationary point is a global minimum, corresponding to a fully synchronized state satisfying $\boldsymbol v_1=\cdots=\boldsymbol v_n$.
Therefore, the central question is:

\begin{quote}
\textit{Which graphs induce benign nonconvexity in the optimization problem \eqref{E_1}?}
\end{quote}
The same problem admits a natural interpretation as a dynamical system. 
To see this, parametrize each unit vector as 
$\boldsymbol v_i=(\cos\theta_i,\sin\theta_i)$. The objective then becomes
\begin{equation}\label{E_2}
\min_{\theta
\in\mathbb{R}^n}
E_G(\boldsymbol\theta)
=
\frac{1}{2}
\sum_{1\le i,j\le n}
A_{ij}\bigl(1-\cos(\theta_i-\theta_j)\bigr).
\end{equation}
The function $E_G(\boldsymbol \theta)$ serves as a potential for the homogeneous Kuramoto model. 
Indeed, since $A$ is symmetric,
\[
\frac{\partial E_G(\boldsymbol \theta)}{\partial \theta_i}
=
\sum_{j=1}^n A_{ij}\sin(\theta_i-\theta_j),
\ \ 1\le i\le n,
\]
and hence its negative gradient flow is
\begin{equation}\label{kuramoto}
    \frac{d\theta_i}{dt}
    =
    -\frac{\partial E_G(\boldsymbol \theta)}{\partial \theta_i}
    =
    -\sum_{j=1}^n A_{ij}\sin(\theta_i-\theta_j),
\ \ 1\le i\le n.
\end{equation}
This is precisely the homogeneous Kuramoto model in a co-rotating frame.
Indeed, the homogeneous Kuramoto model with identical natural frequencies
$\omega_i\equiv\omega$ is
\[
\frac{d\phi_i}{dt}
=
\omega
-
\sum_{j=1}^n A_{ij}\sin(\phi_i-\phi_j),
\qquad 1\le i\le n.
\]
Under the change of variables
$\theta_i=\phi_i-\omega t$, the common drift is removed, yielding
\eqref{kuramoto}.

To relate the benign nonconvexity of \eqref{E_1} to the asymptotic behavior of solutions of \eqref{kuramoto}, we introduce the following definition:
\begin{definition}[Global synchronizing graphs]
    We say a graph $G$ is globally synchronizing if for almost every initial condition, the solution $\theta(t)$ to \eqref{kuramoto} converges to the fully synchronized state $\lim_{t \to \infty} \boldsymbol \theta(t) = \boldsymbol \theta_0$ where all angles are equal modulo \(2\pi\).
\end{definition}
Throughout the paper, when we say that a graph is globally synchronizing, we implicitly assume that it is connected.

In fact, if the function \eqref{E_2} (equivalently \eqref{E_1}) is benign, in the sense that every second-order stationary point is global optimum, then the graph \(G\) is globally synchronizing. Inspired by second-order rigidity in rigidity theory, we call a graph G \emph{second-order globally synchronizing} if the energy function \eqref{E_2} on G is benign. In summary,

\begin{center}
    \textit{$E_G(\boldsymbol \theta)$ is benign, i.e., $G$ is second-order globally synchronizing
    $\Rightarrow$
    $G$ is globally synchronizing.}
\end{center}

This implication follows from standard properties of real-analytic gradient flows. 
The \L{}ojasiewicz gradient inequality~\cite{Lageman2007} guarantees that every trajectory of \eqref{kuramoto} converges to an equilibrium point. 
If an equilibrium is not a second-order stationary point, then its Hessian has a negative eigenvalue, which corresponds to an unstable direction of the gradient flow.
By a standard center-stable manifold argument~\cite{Shub1987}, the set of initial conditions converging to such equilibria has measure zero. 
Consequently, almost every trajectory converges to a second-order stationary point. 
See Lemma~A.1 in~\cite{Geshkovski2025} for a detailed exposition.

The converse is false. For example, the cycle graph $C_4$ is known to be globally synchronizing, yet it admits a second-order stationary point that is not fully synchronized. Namely, the configuration in which the four oscillators form a regular square is an equilibrium with a zero Hessian matrix. Thus, global synchronization does not in general imply that every second-order stationary point is synchronized. Nevertheless, for establishing global synchronization, one may prefer to work with second-order global synchronization, considering it is easier to verify as it requires only first- and second-order optimality conditions, without involving higher-order information.

Beyond its role in nonconvex optimization, this problem also carries physical and biological significance, and more recently, machine learning significance as well.
As a physical phenomenon, the discovery of synchronization dates back to Huygens' seventeenth-century observations of coupled pendulum clocks \cite{PikovskyRosenblumKurths2001}. It also appears throughout nature, for instance in the collective flashing of fireflies in wetlands and in coordinated collective motion, such as the flocking of birds in the sky. Understanding how network topology shapes the long-term dynamics of coupled oscillators, or more generally interacting agents, is a central question in the theory of synchronization and network science.
From an optimization perspective, Ling et al.~\cite{ling2019landscape} present the connection between global synchronization and benign nonconvexity of the rank-two Burer--Monteiro factorization of the semidefinite relaxation for $\mathbb{Z}_2$ synchronization \cite{bandeira2016low}, a problem of recovering a latent signal from noisy or incomplete pairwise measurements.
More recently, connections have been drawn between Kuramoto-type dynamics and dynamical systems arising in transformer architectures. The evolution of token representations across layers can be described by interacting particle systems whose dynamics exhibit alignment and clustering phenomena reminiscent of synchronization. We refer the reader to~\cite{geshkovski2023emergence,Geshkovski2025} for further discussion.

\subsection{Prior work}
Several results are already known characterizing globally synchronizing graphs and second-order globally synchronizing graphs (equivalently, for which \eqref{E_2} exhibits benign nonconvexity). These results can be summarized into the following regimes.
\newline

\noindent
\textit{Regime 1: Minimum-degree proximity to complete graphs}
\newline

\noindent
Recent work ~\cite{taylor2012there,ling2019landscape,lu2021synchronization,kassabov2021sufficiently,canale2022weighted,yoneda2021lower} has focused on characterizing critical thresholds $\mu_c\cdot (n-1)$ for the minimum degree of a graph, above which global synchronization is guaranteed for graphs of arbitrary structure, and below which there exists at least one counterexample.
This can be seen as starting from a complete graph and deleting a fraction of the neighbors of each node, but not too many: as long as more than $75\%$ of the $n-1$ neighbors are retained, global synchrony remains.
The best known bounds on the critical threshold $\mu_c$ currently satisfy $0.6875\leq \mu_c \leq 0.75$ where the upper bound is due to \cite{kassabov2021sufficiently}, and the lower bound follows from \cite{canale2022weighted}.
It has been conjectured that the tight lower bound is $(0.75-\epsilon)(n-1)$ for any $\epsilon>0$ \cite{bandeira2025randomstrasse101openproblems2024}.
Methods beyond linear stability analysis may be needed to establish this bound.
During the preparation of the current version of this manuscript, Lepsveridze and Zhang~\cite{lepsveridze2026graphs} proved this conjecture is false by identifying the four-cluster configuration underlying the $3/4$ barrier and ruling it out as a local minimum through a higher-order perturbation argument.
The exact value of the critical threshold $\mu_c$ remains open; in particular, it is unknown whether the current lower bound $11/16$ is sharp.
\newline

\noindent
\textit{Regime 2: Spectral proximity to complete graphs}
\newline

\noindent
Another regime emerges when random graphs are considered.
The study of global synchronization on random graphs was partially motivated by its connection between the energy landscape \eqref{E_2} and the rank-two Burer--Monteiro factorization of the semidefinite program for $\mathbb Z_2$ synchronization. To see how these two problems are related, represent each variable $\theta_i\in\mathbb R$ by the unit vector
\[
\mathbf q_i=(\cos\theta_i,\sin\theta_i)\in\mathbb R^2,
\ \ 1\le i\le n.
\]
and stack these vectors as the rows of a matrix $\mathbf Q\in\mathbb R^{n\times 2}$, so that
\[
\|\mathbf Q_{i:}\|=1,
\ \ 1\le i\le n.
\]
Under this parametrization, minimizing \eqref{E_2} is equivalent to 
\begin{equation}\label{rank-two-bm}
\max_{\mathbf Q\in\mathbb{R}^{n\times 2}}
\operatorname{Tr}(\mathbf Q^\top \mathbf A \mathbf Q)
\quad
\text{s.t.}\quad
\|\mathbf Q_{i:}\|^2=1,\ \ 1\le i\le n.
\end{equation}

\noindent
Let $\boldsymbol{X} = \mathbf Q\mathbf Q^T$. Then $\boldsymbol{X}$ is positive semidefinite with 
$\operatorname{rank}(\boldsymbol{X}) \le 2$, and the problem~\eqref{rank-two-bm} 
can be equivalently written as
\begin{equation}\label{rank2-sdp}
\max_{\boldsymbol{X}\in \mathbb{R}^{n\times n}} 
\operatorname{Tr}(\mathbf A\boldsymbol{X})
\quad \text{s.t.} \quad \boldsymbol{X} \succeq 0,\;
\operatorname{rank}(\boldsymbol{X}) \le 2,\;
\mathbf X_{ii} = 1,\; 1 \le i \le n.
\end{equation}
Clearly, removing the rank constraint $\operatorname{rank}(\boldsymbol{X}) \le 2$ from \eqref{rank2-sdp}, we obtain semidefinite relaxation of the following combinatorial quadratic form for solving $\mathbb Z_2$ synchronization:
\begin{equation}
    \max_{\boldsymbol{z}\in\{\pm 1\}^n}\boldsymbol{z}^TA\boldsymbol{z}.
\end{equation}
In the context of $\mathbb Z_2$ synchronization, one may consider the noisy observation matrix
\begin{equation}\label{gaussian-matrix}
    \boldsymbol A = \boldsymbol z \boldsymbol z^\top + \sigma \boldsymbol W,
\end{equation}
where $\boldsymbol z\in\{\pm1\}^n$ denotes the ground-truth signal and $\boldsymbol W$ is a Gaussian Wigner matrix~\cite{bandeira2016low}.
Here $\boldsymbol{z}$ corresponds to the ground truth of labeling of nodes where $\mathbf z_i=1$ if it belongs to one group and $\mathbf z_i=-1$ if it belongs to another group. The matrix $\boldsymbol{z}\boldsymbol{z}^T$ then tells for each pair of nodes $i,j$, whether they are assigned into the same group or not. More precisely, if the $(i,j)$ entry of $\boldsymbol{z}\boldsymbol{z}^T$ equals to $1$ then they are assigned into the same group and if it equals to $-1$ then they are assigned into two different groups. The task in the context of $\mathbb Z_2$ synchronization is then to find the maximal strength $\sigma$ of the random perturbation on the pairwise assignment $\boldsymbol{z}\boldsymbol{z}^T$ such that rank-two Burer-Monteiro formalism \eqref{rank-two-bm} still recovers the ground truth, in the sense that every local maximum of \eqref{rank-two-bm} is global\footnote{It is clear that when $\sigma=0$, it successfully recovers the ground truth.}.
Without loss of generality, one can assume $\boldsymbol{z}=\boldsymbol{1}_n$. Then, recovering the ground truth is equivalent to asking how much random perturbation on the edges of the complete graph can be tolerated while preserving global synchrony\footnote{When $\sigma=0$, the adjacency $A=\boldsymbol{1}\boldsymbol{1}^T$ is exactly adjacency matrix of complete graph. It is clear that the diagonal elements $A_{ii}$ do not affect the geometry of landscape since it only adds a constant to the energy function.}.

Substantial progress on this question has been made. In~\cite{bandeira2016low}, it was shown that, under a sufficiently small noise level, every second-order stationary point of the rank-two Burer--Monteiro factorization is rank one and hence globally optimal. This rank-deficiency approach, however, does not yield the optimal noise threshold. More recently, Rakoto Endor and Waldspurger~\cite{endor2024benign} took a different perspective and established a sharp deterministic condition, in terms of the condition number of an associated Laplacian matrix, that guarantees a benign landscape. Building on this approach, McRae~\cite{mcrae2025benign} introduced a diagonally preconditioned Laplacian condition and obtained the asymptotically tight noise bound
\[
    \sigma \leq \sqrt{\frac{n}{(2+\epsilon)\log n}},
\]
for any $\epsilon>0$, which asymptotically matches the information-theoretic threshold for exact recovery. For several other random graph models, including Erd\H{o}s--R\'enyi graphs, random $d$-regular graphs, and random graph processes, and several deterministic graphs such as expander graphs \cite{abdalla2026expander}, global synchronization in suitable regimes can also be established via a condition-number criterion~\cite{mcrae2025benign}.

Taken together, these results show that a graph is second-order globally synchronizing whenever the condition number of its preconditioned Laplacian is sufficiently close to one, which is the condition number of the preconditioned Laplacian of complete graphs. Thus, these results can be viewed as extending second-order global synchronization from complete graphs to graphs that are sufficiently close to complete graphs in this spectral sense.
\newline

\noindent
\textit{A unified view of regimes 1 and 2}
\newline

\noindent
To summarize previous results, in fact, both \textit{Regime 1} and \textit{Regime 2} share the same underlying philosophy: complete graphs globally synchronize, and graphs that are sufficiently close to complete in certain sense also globally synchronize.
\[
\{\textit{Graphs } \approx K_n\}\ \textit{globally synchronize.}
\]
The two regimes differ in how proximity to the complete graph is quantified: Regime 1 uses a combinatorial notion based on minimum degree, whereas Regime 2 uses a spectral notion based on the Laplacian.
\newline

\noindent
\textit{Our contribution: a new structural mechanism}
\newline

\noindent
In this work, we identify an entirely different mechanism for second-order global synchronization, one that does not rely on the close-to-complete regime underlying previous results.
We prove that the comparability graph associated with the partial order induced by a rooted tree is second-order globally synchronizing, and hence globally synchronizing. These graphs are precisely the connected quasi-threshold graphs.

\begin{theorem}\label{main1}
Every connected quasi-threshold graph is second-order globally synchronizing.
\end{theorem}
This implies every connected quasi-threshold graph is globally synchronizing.
Moreover, quasi-threshold graphs can also be defined recursively starting from $K_1$ by repeatedly applying either of the following two operations: taking the disjoint union of existing quasi-threshold graphs, or adding a universal vertex adjacent to every vertex of an existing quasi-threshold graph.

To the best of our knowledge, the result most closely related to ours is our previous result on threshold graphs.
However, the present work reveals that the phenomenon is not specific to threshold graphs, but arises from the considerably more general rooted-tree structure underlying quasi-threshold graphs.
Moreover, the proof technique developed in~\cite{wu2025threshold} does not extend to the results in the current work, as we explain in Section~\ref{sec:beyond-threshold}.











In a nutshell, our proof uncovers a structural mechanism for global synchronization that is not captured by existing criteria based on minimum degree or the condition number of the graph Laplacian. The key is an underlying rooted-tree skeleton: synchronization is first forced within suitable local substructures and then propagates along the skeleton, eventually forcing every second-order stationary point to be fully synchronized.


\subsection{Structure of the paper}
The preliminaries on the Kuramoto model and the vector geometry of equilibria are presented in Section~\ref{sec:pre}.
In Section~\ref{sec:skeleton}, we introduce partially ordered sets and explain how they define comparability graphs.
In Section~\ref{sec:beyond-threshold}, we review previous results on the global synchronization of threshold graphs, and explain why the previous proof strategy fails for quasi-threshold graphs.
In the next section, we introduce geometric twins and the scenarios in which they arise.
In Section~\ref{sec:quasi-threshold}, we present a propagation argument on rooted trees, which shows that quasi-threshold graphs are globally synchronizing.
\subsection{Notation}\label{notation}

\subsubsection{Basic graph notation}
For an undirected graph $G=(V,E)$, a pair $uv\in E$ is called an \emph{edge}. The open neighbourhood of node $u$ is defined as $N(u)=\{v\in V: uv\in E\}$ and the closed neighbourhood of $u$ is defined as
$N[u]=N(u)\cup\{u\}$.
For a directed graph $G=(V,E)$, an ordered pair
$(u,v)\in E$ is called an \emph{arc}.
In an undirected graph $G=(V,E)$, nodes $u$ and $v$ are called \emph{structural closed twins} if $N[u]=N[v]$, and \emph{structural open twins} if $N(u)=N(v)$.
We use the adjective ``structural'' to distinguish these notions from the geometric twins introduced in Section~\ref{sec:geo-twins}. Let $G[B]$ denote the subgraph of $G$ induced by the vertex set $B$.
For a rooted tree $T$ and a node $v\in V(T)$, let $\operatorname{Anc}(v)$ and $\operatorname{Desc}(v)$ denote, respectively, the sets of proper ancestors and proper descendants of $v$, and let $\operatorname{Children}(v)$ denote the set of children of $v$.

\subsubsection{Phase and phasor notation}
For a phase configuration $\boldsymbol\theta = (\theta_1,\ldots,\theta_n) \in \mathbb{R}^n$, we associate to each node $i \in V$ the unit vector $\boldsymbol v_i = (\cos \theta_i, \sin \theta_i) \in \mathbb{S}^1$, which we refer to as the phasor of node $i$. We often identify a phase configuration $\boldsymbol\theta$ with its corresponding phasor configuration $(\boldsymbol v_1,\ldots,\boldsymbol v_n)$. Two nodes $i,j \in V$ are said to synchronize if $\theta_i = \theta_j \pmod{2\pi}$, equivalently, if $\boldsymbol v_i = \boldsymbol v_j$.
For nonzero vectors $\boldsymbol x,\boldsymbol y\in\mathbb{R}^2$, we write $\boldsymbol x \Uparrow \boldsymbol y$ if $\boldsymbol x$ and $\boldsymbol y$ point in the same direction.


\subsubsection{Graph operations}\label{sec:operations}
We summarize the graph operations on two graphs $G=(V_G,E_G)$ and $H=(V_H,E_H)$ used in this paper as follows. Here, $V$ and $E$ denote the vertex and edge sets of the resulting graph.

\begin{center}
\begin{tabular}{lll}
\hline
Operation & Notation & Description \\
\hline
Disjoint union & $G \oplus H$ 
& $V=V_G\cup V_H,\ E=E_G\cup E_H$ \\

Add isolated vertex & $G\,\Iso\, i$ 
& $V=V_G\cup\{i\},\ E=E_G$ \\

Add universal vertex & $G\,\Uni\, i$ 
& $V=V_G\cup\{i\},\ E=E_G\cup\{ui: u\in V_G\}$ \\
\hline
\end{tabular}
\end{center}

\section{Preliminaries on Kuramoto model and vector geometry of equilibria}\label{sec:pre}
In this section, we present the basics of the Kuramoto model, including the first- and second-order stationary conditions of the associated nonconvex optimization landscape, together with their geometric interpretations in the plane. We briefly summarize the basic geometric lemmas here.

\subsection{First-order stationary condition}
Let $G=(V,E)$ be an undirected graph with adjacency matrix $\boldsymbol A$.
We consider the Kuramoto energy
\[
 E_G(\boldsymbol\theta)
=
\frac12 \sum_{i=1}^n \sum_{j=1}^n  A_{ij}\bigl(1-\cos(\theta_i-\theta_j)\bigr).
\]
A straightforward computation yields
\[
\bigl(\nabla E_G(\boldsymbol\theta)\bigr)_j
=
\sum_{i=1}^n A_{ij}\sin(\theta_j-\theta_i),
\qquad 1 \le j \le n .
\]
Therefore, $\boldsymbol \theta$ is a first-order stationary point of $ E_G(\boldsymbol\theta)$ if and only if
\[
\sum_{i=1}^n A_{ij}\sin(\theta_j-\theta_i)=0,
\qquad \forall\, 1 \le j \le n.
\]
The following geometric characterization of first-order stationary points was established in \cite{monzon2005global}. A proof is also given in Lemma 4.1 of \cite{wu2025threshold}.
\begin{lemma}\label{eq-vec}
A state \( \boldsymbol \theta \) is an equilibrium of the Kuramoto model~\eqref{kuramoto} on $G=(V,E)$ with adjacency $A$ {if and only if}, for each \( 1 \le i \le n\),
\begin{equation}\label{formula:eq-ec}
\sum_{j\in N(i)} \boldsymbol v_j = \mu_i \boldsymbol v_i\quad\text{where}\quad \mu_i=\sum_{j\in N(i)}\cos(\theta_j-\theta_i),
\end{equation}
where $\boldsymbol v_i=(\cos\theta_i,\sin\theta_i)$.
\end{lemma}

\subsection{Second-order stationary condition}
The Hessian of $E_G(\boldsymbol{\theta})$ at $\boldsymbol{\theta}$ has entries, for $1 \le i,j \le n$,
\[
(\nabla^2 E_G(\boldsymbol{\theta}))_{i,j}
=
\begin{cases}
-A_{ij}\cos(\theta_i-\theta_j), & i\neq j,\\[1mm]
\displaystyle\sum_{k\in V\setminus\{i\}}
A_{ik}\cos(\theta_i-\theta_k), & i=j.
\end{cases}
\]
If $\boldsymbol\theta$ is a second-order stationary point of $E_G(\boldsymbol\theta)$, then the Hessian is
positive semidefinite. This implies the following lemma.

\begin{lemma}\label{stable-eq}
If $\boldsymbol{\theta}$ is a second-order stationary point of~\eqref{E_2},
then \eqref{formula:eq-ec} holds with $\mu_i\ge 0$ for every $1\le i\le n$.
\end{lemma}

The proof of Lemma~\ref{stable-eq} can be found in Lemma~4.4 of~\cite{wu2025threshold}.

\section{Graph skeletons: comparability graphs of partially ordered sets}\label{sec:skeleton}

\noindent


\noindent
\textbf{Rooted forest, induced graph, transitive closure.}
A \textit{rooted tree} $T$ is a directed graph obtained by
orienting the edges of a tree toward a distinguished vertex
$r$, called the \textit{root}
(shown in black in Figure~\ref{rootedtree}).
Thus, every vertex has a unique directed path to $r$. A \textit{rooted forest} $F$ is a graph whose connected components are rooted trees.
We say a rooted tree $T$ \textit{induces} an undirected graph $G(T)=(V,E)$, if $uv$ is an edge in $E$, if and only if, $u\neq v$ and there exists a directed path from $u$ to $v$ or from $v$ to $u$.
The transitive closure of a rooted tree $T$ is a directed graph $C(T)$ where $u,v$ ($u\neq v$) is an arc if and only if there exists a directed path from $u$ to $v$.
\newline

\noindent
\textbf{Posets and comparability graphs.}
A partially ordered set (poset) is a pair $P=(V,\preceq)$, where $\preceq$ is a reflexive, antisymmetric, and transitive binary relation on $V$. Two distinct elements $u,v\in V$ are said to be comparable if either $u\preceq v$ or $v\preceq u$. The comparability graph of $P$ is the undirected graph on vertex set $V$ in which two distinct vertices are adjacent if and only if they are comparable.

A rooted tree $T$ defines a partial order $P=(V,\preceq)$ by setting
$u\preceq v$ if $u=v$ or there exists a directed path from $u$ to $v$ in $T$.
Then the transitive closure $C(T)$ contains an arc $(u,v)$ if and only if
$u \preceq v$ and $u\neq v$. Moreover, $T$ is the Hasse diagram of $P$, and $G(T)$ is the comparability graph of $P$.

\begin{figure}[H]
\centering
\begin{tikzpicture}[scale=1,
    every node/.style={circle, draw, minimum size=3mm, inner sep=0pt},
    level distance=5mm,
    treeedge/.style={black, line width=1pt, ->, >=stealth, shorten >=2pt},
    uedge/.style={black, line width=1pt}]

\begin{scope}[xshift=0cm]
\node[fill=black, label=above:{6}] (r3) at (0,0) {};
\node[fill=white, label=above:{4}] (a3) at (-1.2,-1.2) {};
\node[fill=white, label=above:{5}] (b3) at (1.2,-1.2) {};
\node[fill=white, label=below:{1}] (c3) at (-2.0,-2.4) {};
\node[fill=white, label=below:{2}] (d3) at (-0.4,-2.4) {};
\node[fill=white, label=below:{3}] (e3) at (1.2,-2.4) {};

\draw[treeedge] (a3) -- (r3);
\draw[treeedge] (b3) -- (r3);
\draw[treeedge] (c3) -- (a3);
\draw[treeedge] (d3) -- (a3);
\draw[treeedge] (e3) -- (b3);

\node[draw=none, rectangle] at (0,-3.3) {(A)};
\end{scope}


\begin{scope}[xshift=5.5cm]

\node[fill=white, label=above:{6}] (r2) at (0,0) {};
\node[fill=white, label=above:{4}] (a2) at (-1.2,-1.2) {};
\node[fill=white, label=above:{5}] (b2) at (1.2,-1.2) {};
\node[fill=white, label=below:{1}] (c2) at (-2.0,-2.4) {};
\node[fill=white, label=below:{2}] (d2) at (-0.4,-2.4) {};
\node[fill=white, label=below:{3}] (e2) at (1.2,-2.4) {};

\draw[uedge] (r2)--(a2);
\draw[uedge] (r2)--(b2);
\draw[uedge] (a2)--(c2);
\draw[uedge] (a2)--(d2);
\draw[uedge] (b2)--(e2);

\draw[line width=1.2pt, bend right=40] (r2) to (c2);
\draw[line width=1.2pt, bend right=15] (r2) to (d2);
\draw[line width=1.2pt, bend right=10] (r2) to (e2);

\node[draw=none, rectangle] at (0,-3.3) {(B)};

\end{scope}

\begin{scope}[xshift=11cm]
\node[fill=white, label=above:{6}] (r2) at (0,0) {};
\node[fill=white, label=above:{4}] (a2) at (-1.2,-1.2) {};
\node[fill=white, label=above:{5}] (b2) at (1.2,-1.2) {};
\node[fill=white, label=below:{1}] (c2) at (-2.0,-2.4) {};
\node[fill=white, label=below:{2}] (d2) at (-0.4,-2.4) {};
\node[fill=white, label=below:{3}] (e2) at (1.2,-2.4) {};

\draw[treeedge] (a2) -- (r2);
\draw[treeedge] (b2) -- (r2);
\draw[treeedge] (c2) -- (a2);
\draw[treeedge] (d2) -- (a2);
\draw[treeedge] (e2) -- (b2);

\draw[treeedge, bend left=40] (c2) to (r2);
\draw[treeedge, bend left=15] (d2) to (r2);
\draw[treeedge, bend left=10] (e2) to (r2);

\node[draw=none, rectangle] at (0,-3.3) {(C)};

\end{scope}

\end{tikzpicture}
\caption{(A) Rooted tree $T$; (B) Induced graph $G(T)$ of $T$; (C) Transitive closure $C(T)$ of $T$.}
\label{rootedtree}
\end{figure}

\noindent
\textbf{Quasi-threshold graphs as comparability graphs.}
The following characterization connects quasi-threshold graphs with
the rooted-tree structure introduced above.
\begin{theorem}[\cite{yan1996quasi}, Theorem 3]
A graph is quasi-threshold if and only if it is the comparability graph
of the partial order induced by a rooted forest. Moreover, it is connected
if and only if the underlying rooted forest is a rooted tree. \end{theorem}
Quasi-threshold graphs also admit an equivalent recursive characterization: starting from a single vertex, they can be generated by repeated disjoint unions and by adding a universal vertex. This operational viewpoint will occasionally be useful for comparison with threshold graphs, while our main arguments will rely on the rooted-tree representation above.

For example, the graph in Figure~\ref{rootedtree}(B) has the following construction sequence, where the graph-operation notation is defined in Section~\ref{sec:operations}.
\begin{figure}[H]
\centering
\begin{tikzpicture}[scale=1,
    every node/.style={circle, draw, minimum size=3mm, inner sep=0pt},
    level distance=5mm,
    treeedge/.style={black, line width=1pt, ->, >=stealth, shorten >=2pt},
    uedge/.style={black, line width=1pt}]

\begin{scope}[xshift=0cm]
\node[draw=none, rectangle] at (0,-1.3) {$\left(\left(\left(1  \oplus  2 \right) \; \Uni \; 4\right)\oplus \left( 3 \;\Uni\; 5 \right) \right) \;\Uni\; 6$};
\end{scope}
\end{tikzpicture}
\end{figure}

A rooted tree is called a \emph{rooted caterpillar} if
every node has at most one child that is not a leaf.
Equivalently, all nodes with children lie on a single
directed path ending at the root.
Here, a leaf is a node with no children.

\begin{theorem}
A graph is a connected threshold graph if and only if it is the comparability graph of a poset defined by a rooted caterpillar.
\end{theorem}

\begin{figure}[H]
\centering
\begin{tikzpicture}[scale=1,
    every node/.style={circle, draw, minimum size=3mm, inner sep=0pt},
    level distance=5mm,
    treeedge/.style={black, line width=1pt, ->, >=stealth, shorten >=2pt},
    uedge/.style={black, line width=1pt}]

\begin{scope}[xshift=0cm]
\node[fill=black, label=above:{4}] (r3) at (0,0) {};
\node[fill=white, label=above:{3}] (a3) at (-1.2,-1.2) {};
\node[fill=white, label=below:{1}] (c3) at (-2.0,-2.4) {};
\node[fill=white, label=below:{2}] (d3) at (-0.4,-2.4) {};

\draw[treeedge] (a3) -- (r3);
\draw[treeedge] (c3) -- (a3);
 \draw[treeedge] (d3) -- (a3);

\node[draw=none, rectangle] at (-1.2,-3.3) {(A)};
\end{scope}

\begin{scope}[xshift=5cm]
\node[draw=none, rectangle] at (0,-1.3) {$\left(\left(1 \; \Iso \; 2 \right)\; \Uni \; 3\right)\; \Uni \; 4$};
\node[draw=none, rectangle] at (0,-3.3) {(B)};
\end{scope}
\begin{scope}[xshift=12cm]

\node[fill=white, label=above:{4}] (r2) at (0,0) {};
\node[fill=white, label=above:{3}] (a2) at (-1.2,-1.2) {};
\node[fill=white, label=below:{1}] (c2) at (-2.0,-2.4) {};
\node[fill=white, label=below:{2}] (d2) at (-0.4,-2.4) {};

\draw[uedge] (r2)--(a2);
\draw[uedge] (a2)--(c2);
\draw[uedge] (a2)--(d2);

\draw[line width=1.2pt, bend right=40] (r2) to (c2);
\draw[line width=1.2pt, bend right=15] (r2) to (d2);

\node[draw=none, rectangle] at (-1.2,-3.3) {(C)};

\end{scope}
\end{tikzpicture}
\caption{This figure illustrates why the threshold graph in (C), constructed via (B), is a comparability graph of the poset defined by the skeleton in (A). Starting from node $1$, we add node $2$ as an isolated vertex, obtaining the subgraph represented by the innermost parentheses. We then add nodes $3$ and $4$ as universal vertices, thereby obtaining the graph in (C).}
\end{figure}

\section{Beyond Threshold graphs: Proof difficulty and New Ideas}\label{sec:beyond-threshold}
In this section, we explain why the proof of global synchronization for threshold graphs developed in \cite{wu2025threshold} does not extend directly to quasi-threshold graphs.
We begin by reinterpreting the argument in \cite{wu2025threshold} from the perspective of comparability graphs of posets induced by rooted caterpillars.
Consider, for example, the threshold graph constructed
from an initial vertex $A$ using the sequence
\texttt{10101011001}.
The bits encode the eleven subsequent vertex additions:
$1$ denotes a universal vertex and $0$ an isolated vertex.
These added vertices are labeled by uppercase and lowercase
letters, respectively, as shown below.
\begin{center}
\begin{tikzpicture}[>=stealth,thick]
  \foreach \i/\digit/\letter in {
    0//A,
    1/\texttt{1}/B,
    2/\texttt{0}/c,
    3/\texttt{1}/D,
    4/\texttt{0}/e,
    5/\texttt{1}/F,
    6/\texttt{0}/g,
    7/\texttt{1}/H,
    8/\texttt1/I,
    9/\texttt{0}/j,
    10/\texttt{0}/k,
    11/\texttt{1}/L
  }{
    \node (num\i) at (0.6*\i,0) {\digit}; 
    \node (let\i) at (0.6*\i,-1) {\letter};
    \draw[->] (num\i.south) -- (let\i.north);
  }
\end{tikzpicture}
\end{center}

To construct its underlying skeleton, we connect all universal vertices (including the initial vertex if the first bit is $\texttt{1}$) as a path where later-added universal vertices are closer to the root than those added earlier. Then, attach each isolated vertex to the nearest subsequent universal vertex in the construction sequence. Specify the last added universal node $L$ as the root of the caterpillar. The resulting rooted caterpillar is shown on the left-hand side of Figure \ref{fig:rooted-caterpillar-example}.

By connecting every node to all of its ancestors in the caterpillar, we recover the threshold graph shown on the right-hand side of Figure~\ref{fig:rooted-caterpillar-example}.

\begin{figure}[htbp]
\centering
\begin{tikzpicture}[
    x=1cm,y=1cm,font=\small,
    cat vertex/.style={circle,draw=black,fill=white,
        line width=0.4pt,minimum size=2.3mm,
        inner sep=0pt,outer sep=0pt},
    cat tree/.style={draw=black,line width=0.65pt},
    cat arrow/.style={cat tree,
        -{Stealth[length=1.45mm,width=1.05mm]}},
    cat closure/.style={draw=black!100,line width=0.44pt},
    every node/.style={inner sep=1pt}
]
    \foreach \panel/\xx in {1/0,2/7.50} {
        \begin{scope}[xshift=\xx cm]
            \node[cat vertex] (A) at (0.00,0.00) {};
            \node[cat vertex] (B) at (0.82,0.76) {};
            \node[cat vertex] (D) at (1.64,1.52) {};
            \node[cat vertex] (F) at (2.46,2.28) {};
            \node[cat vertex] (H) at (3.28,3.04) {};
            \node[cat vertex] (I) at (4.10,3.80) {};
            \node[cat vertex] (L) at (4.92,4.56) {};
            \node[cat vertex] (c) at (1.64,0.76) {};
            \node[cat vertex] (e) at (2.46,1.52) {};
            \node[cat vertex] (g) at (3.28,2.28) {};
            \node[cat vertex] (j) at (4.92,3.80) {};
            \node[cat vertex] (k) at (5.58,3.80) {};

            \ifnum\panel=1\relax
                \node[cat vertex,fill=black] at (L) {};
            \fi

            \begin{pgfonlayer}{background}
                \ifnum\panel=2\relax
                    \foreach \u/\v/\bend in {
                        D/A/60,
                        F/A/60,
                        H/A/60,
                        I/A/60,
                        L/A/60,
                        F/B/60,
                        H/B/60,
                        I/B/60,
                        L/B/60,
                        H/D/60,
                        I/D/60,
                        L/D/60,
                        I/F/60,
                        L/F/60,
                        L/H/60,
                        F/c/-70,
                        H/c/-70,
                        I/c/-70,
                        L/c/-70,
                        H/e/-70,
                        I/e/-70,
                        L/e/-70,
                        I/g/-70,
                        L/g/-70%
                    } {
                        \draw[cat closure]
                            let \p1=(\u.center), \p2=(\v.center),
                                \n1={atan2(\y2-\y1,\x2-\x1)
                                    +ifthenelse(\bend>0,-90,90)-\bend},
                                \n2={veclen(\x2-\x1,\y2-\y1)
                                    /(2*sin(abs(\bend)))}
                            in (\u.center)
                                arc[start angle=\n1,
                                    delta angle={2*\bend},radius=\n2];
                    }
                \fi

                \foreach \parent/\child in {
                    B/A,
                    D/B,
                    F/D,
                    H/F,
                    I/H,
                    L/I,
                    D/c,
                    F/e,
                    H/g,
                    L/j,
                    L/k%
                } {
                    \ifnum\panel=1\relax
                        \draw[cat arrow] (\child) -- (\parent);
                    \else
                        \draw[cat tree] (\parent) -- (\child);
                    \fi
                }
            \end{pgfonlayer}

            \foreach \v in {A,B,D,F,H,I,c,e,g,j} {
                \node[anchor=east,fill=white,inner sep=0.6pt]
                    at ($(\v.west)+(-0.09,0)$) {$\v$};
            }
            \node[right=3pt] at (k.east) {$k$};
            \node[above=3pt] at (L.north) {$L$};

            \ifnum\panel=1\relax
            \fi
        \end{scope}
    }
\end{tikzpicture}
\caption{A rooted caterpillar (left) and its corresponding comparability graph (right).}
\label{fig:rooted-caterpillar-example}
\end{figure}

Recall that the proof of global synchronization for threshold graphs in \cite{wu2025threshold} is inductive.
On the present example, the argument proceeds by propagating synchronization along the spine of the graph, starting from the root $L$. 
Specifically, the first step is to establish synchronization within the group $\{L, j, k\}$.
Next, one shows that the previously added universal node $I$ synchronizes with $L$. 
Continuing in this manner, synchronization is propagated iteratively down the spine, until all nodes are synchronized.
In other words, the process proceeds from the highest node on the spine (the root $L$) to the lowest node $A$, establishing synchronization step by step along the way.

The key to propagation from the root down the spine is the following.
Consider two consecutive universal vertices along the spine.
If no isolated vertices were inserted between their additions,
they are structural closed twins and synchronize at any second-order
stationary point by Corollary~5.4 of \cite{wu2025threshold}.
Otherwise, their closed neighborhoods differ precisely by the
intervening isolated vertices.
At the relevant induction step, these extra vertices have already
been shown to synchronize with the later-added universal vertex.
They therefore form a synchronous pendant extension, and
Lemma~5.9 of \cite{wu2025threshold} implies that the two universal
vertices also synchronize.
This allows synchronization to propagate down the spine.

The obstruction is the branching structure of a general rooted tree: the neighborhood difference between a parent and a child may contain entire sibling subtrees, so the synchronous pendant-extension argument from the threshold-graph case no longer applies directly.
However, this approach appears to break down already at the first step, as it would require showing that, at any second-order stationary point, the root node $A$ synchronizes with its children $B$ and $C$. 
At present, it is not clear whether the difference between $N[A]$ and $N[B]$ forms a synchronous pendant extension, as in the threshold graph setting.

\begin{figure}[htbp]
\centering
\begin{tikzpicture}[
    x=1cm,y=1cm,font=\small,
    rc vertex/.style={circle,draw=black,fill=white,
        line width=0.40pt,minimum size=2.3mm,
        inner sep=0pt,outer sep=0pt},
    rc tree/.style={draw=black,line width=0.65pt},
    rc arrow/.style={rc tree,
        -{Stealth[length=1.45mm,width=1.05mm]}},
    rc closure/.style={draw=black!100,line width=0.44pt},
    every node/.style={inner sep=1pt}
]
    \foreach \panel/\xx in {1/0,2/7.85} {
        \begin{scope}[xshift=\xx cm]
            \node[rc vertex] (A) at (0.00,0.00) {};
            \node[rc vertex] (B) at (-1.15,-0.82) {};
            \node[rc vertex] (M) at (0.00,-0.82) {};
            \node[rc vertex] (C) at (1.45,-0.82) {};
            \node[rc vertex] (b1) at (-1.97,-1.64) {};
            \node[rc vertex] (b2) at (-2.72,-2.46) {};
            \node[rc vertex] (b3) at (-1.97,-2.46) {};
            \node[rc vertex] (b4) at (-2.72,-3.28) {};
            \node[rc vertex] (m1) at (-0.63,-1.64) {};
            \node[rc vertex] (m2) at (0.12,-1.64) {};
            \node[rc vertex] (m3) at (-0.63,-2.46) {};
            \node[rc vertex] (m4) at (-1.22,-3.28) {};
            \node[rc vertex] (m5) at (-0.29,-3.28) {};
            \node[rc vertex] (c1) at (0.87,-1.64) {};
            \node[rc vertex] (c2) at (2.32,-1.64) {};
            \node[rc vertex] (c3) at (0.87,-2.46) {};
            \node[rc vertex] (c4) at (0.28,-3.28) {};
            \node[rc vertex] (c5) at (1.45,-3.28) {};
            \node[rc vertex] (c6) at (2.32,-2.46) {};
            \node[rc vertex] (c7) at (3.08,-2.46) {};

            \ifnum\panel=1\relax
                \node[rc vertex,fill=black] at (A) {};
            \fi

            \begin{pgfonlayer}{background}
                \ifnum\panel=2\relax
                    \foreach \u/\v/\bend in {
                        A/b1/42.00,
                        B/b2/30.00,
                        A/b2/42.00,
                        B/b3/24.00,
                        A/b3/66.00,
                        b1/b4/18.00,
                        B/b4/6.00,
                        A/b4/66.00,
                        A/m1/42.00,
                        A/m2/-42.00,
                        M/m3/12.00,
                        A/m3/42.00,
                        m1/m4/30.00,
                        M/m4/36.00,
                        A/m4/42.00,
                        m1/m5/-30.00,
                        M/m5/24.00,
                        A/m5/54.00,
                        A/c1/-42.00,
                        A/c2/-42.00,
                        C/c3/12.00,
                        A/c3/-48.00,
                        c1/c4/30.00,
                        C/c4/36.00,
                        A/c4/-36.00,
                        c1/c5/-30.00,
                        C/c5/-30.00,
                        A/c5/-42.00,
                        C/c6/-24.00,
                        A/c6/-18.00,
                        C/c7/-30.00,
                        A/c7/-54.00%
                    } {
                        \draw[rc closure]
                            let \p1=(\u.center), \p2=(\v.center),
                                \n1={atan2(\y2-\y1,\x2-\x1)
                                    +ifthenelse(\bend>0,-90,90)-\bend},
                                \n2={veclen(\x2-\x1,\y2-\y1)
                                    /(2*sin(abs(\bend)))}
                            in (\u.center)
                                arc[start angle=\n1,
                                    delta angle={2*\bend},radius=\n2];
                    }
                \fi

                \foreach \parent/\child in {
                    A/B,
                    A/M,
                    A/C,
                    B/b1,
                    b1/b2,
                    b1/b3,
                    b2/b4,
                    M/m1,
                    M/m2,
                    m1/m3,
                    m3/m4,
                    m3/m5,
                    C/c1,
                    C/c2,
                    c1/c3,
                    c3/c4,
                    c3/c5,
                    c2/c6,
                    c2/c7%
                } {
                    \ifnum\panel=1\relax
                        \draw[rc arrow] (\child) -- (\parent);
                    \else
                        \draw[rc tree] (\parent) -- (\child);
                    \fi
                }
            \end{pgfonlayer}

            \node[above=3pt] at (A.north) {$A$};
            \node[above left=2pt,fill=white] at (B.north west) {$B$};
            \node[above right=2pt,fill=white] at (C.north east) {$C$};

            \ifnum\panel=1\relax
            \fi
        \end{scope}
    }
\end{tikzpicture}
\caption{A rooted tree (left) and its corresponding comparability graph (right).}
\label{fig:rooted-tree-comparability}
\end{figure}

The propagation direction can instead be reversed: leaves sharing the same parent are forced to synchronize first, and this local synchronization then propagates upward layer by layer. 
In this way, at any second-order stationary point, synchronization ultimately spreads to the entire graph.
The detailed proof is given in Section \ref{sec:quasi-threshold}, while Section \ref{sec:geo-twins} first establishes the lemmas needed for the proof.

%% file: 02_quasi.tex
\section{Phasor geometry of geometric twins}\label{sec:geo-twins}
\subsection{Two basic geometric facts}
We begin with a geometric lemma that characterizes all
possible relative positions of two unit vectors satisfying a pair of linear relations.

\begin{lemma}\label{geo}
Let $\mathbf{v}_a, \mathbf{v}_b \in \mathbb{S}^1$. Suppose there exist a vector $\boldsymbol q \in \mathbb{R}^2$ and scalars $\mu_a, \mu_b \in \mathbb{R}$ such that
\begin{equation}\label{acase1}
  \mathbf{v}_b + \boldsymbol q = \mu_a \mathbf{v}_a \qquad \mathrm{and} \qquad \mathbf{v}_a + \boldsymbol q = \mu_b \mathbf{v}_b.
\end{equation}
Then the positions of $\mathbf{v}_a$ and $\mathbf{v}_b$ fall into one of the following three cases:
\begin{enumerate}
    \item $\mathbf{v}_a = \mathbf{v}_b$ and $\mu_a = \mu_b$;
    \item $\mathbf{v}_a = -\mathbf{v}_b$, $\mu_a + \mu_b = -2$, and $(\mu_a, \mu_b) \neq (-1,-1)$;
    \item $\mathbf{v}_a, \mathbf{v}_b \in \mathbb{S}^1$, $\mu_a = \mu_b = -1$, and $\mathbf{v}_a + \mathbf{v}_b + \boldsymbol q = 0$.
\end{enumerate}
Alternatively, suppose that
\begin{equation}\label{acase2}
   \boldsymbol q  =  \mu_a \mathbf{v}_a \qquad \mathrm{and} \qquad \boldsymbol q  = \mu_b \mathbf{v}_b.
\end{equation}
Then the positions of $\mathbf{v}_a$ and $\mathbf{v}_b$ fall into one of the following three cases:
\begin{enumerate}
    \item $\mathbf{v}_a = \mathbf{v}_b$ and $\mu_a = \mu_b \ne 0$;
    \item $\mathbf{v}_a = -\mathbf{v}_b$ and $\mu_a = -\mu_b \ne 0$;
    \item $\mathbf{v}_a, \mathbf{v}_b \in \mathbb{S}^1$ and $\mu_a = \mu_b = 0$.
\end{enumerate}
\end{lemma}
The corresponding feasible regions for $(\mu_a,\mu_b)$ are illustrated in Figure~\ref{fig:relations}.
\begin{figure}[H]
\centering

\begin{minipage}[c]{0.45\textwidth}
\centering
\begin{tikzpicture}[scale=0.65]
    \draw[thick, ->] (-3, 0) -- (3, 0) node[right] {$\mu_a$};
    \draw[thick, ->] (0, -3) -- (0, 3) node[above] {$\mu_b$};
    \draw[thick, black] (-3, -3) -- (3, 3);
    \draw[thick, black] (-3,1) -- (1,-3);
    \draw[fill=white, draw=black, thick] (-1,-1) circle (3pt);
\end{tikzpicture}
\end{minipage}
\hfill
\begin{minipage}[c]{0.45\textwidth}
\centering
\begin{tikzpicture}[scale=0.65]
    \draw[thick, ->] (-3, 0) -- (3, 0) node[right] {$\mu_a$};
    \draw[thick, ->] (0, -3) -- (0, 3) node[above] {$\mu_b$};
    \draw[thick, black] (-3, -3) -- (3, 3);
    \draw[thick, black] (-3, 3) -- (3, -3);
    \draw[fill=white, draw=black, thick] (0,0) circle (3pt);
\end{tikzpicture}
\end{minipage}

\caption{Feasible regions (the lines) of $(\mu_a,\mu_b)$ corresponding to cases~\eqref{acase1} and~\eqref{acase2}, respectively.}
\label{fig:relations}
\end{figure}
Here, we prove the result for the case \eqref{acase2}. The proof for the case \eqref{acase1} can be found in Lemma~5.2 of \cite{wu2025threshold}.

\begin{proof}
From \eqref{acase2} and
$\|\mathbf v_a\|=\|\mathbf v_b\|=1$, we obtain
\[
    \|\boldsymbol q\|=|\mu_a|=|\mu_b|.
\]
Thus, $\mu_a$ and $\mu_b$ are either both zero or both nonzero.
The following two cases therefore exhaust all possibilities.

\textit{Case 1.} If $\mu_a=\mu_b=0$, then
$\boldsymbol q=\mathbf 0$.

\textit{Case 2.} If $\mu_a,\mu_b\ne 0$, then
\[
    \mathbf v_a
    =\frac{1}{\mu_a}\boldsymbol q
    =\frac{\mu_b}{\mu_a}\mathbf v_b.
\]
Since $|\mu_a|=|\mu_b|$, we distinguish two subcases.
If $\mu_a=\mu_b$, then $\mathbf v_a=\mathbf v_b$.
If $\mu_a=-\mu_b$, then $\mathbf v_a=-\mathbf v_b$,
and one of $\mu_a,\mu_b$ is negative.
\end{proof}

The two geometric configurations characterized in Lemma~\ref{geo} will arise repeatedly in our synchronization arguments, so we introduce the following terminology.
\begin{definition}[Geometrically open and closed twins]
At a state $\boldsymbol{\theta}$ of the Kuramoto model on
$G$, let $i,j\in V(G)$ be distinct nodes, and let
$\mathbf q\in\mathbb R^2$ be a prescribed interaction vector.

\begin{enumerate}
    \item The nodes $i$ and $j$ are called
    \emph{geometrically closed twins relative to $\mathbf q$}
    if there exist $\alpha_i,\alpha_j\in\mathbb R$ such that
    \[
        \mathbf v_j+\mathbf q=\alpha_i\mathbf v_i,
        \qquad
        \mathbf v_i+\mathbf q=\alpha_j\mathbf v_j.
    \]

    \item The nodes $i$ and $j$ are called
    \emph{geometrically open twins relative to $\mathbf q$}
    if there exist $\alpha_i,\alpha_j\in\mathbb R$ such that
    \[
        \mathbf q=\alpha_i\mathbf v_i=\alpha_j\mathbf v_j.
    \]
\end{enumerate}

In either case, we add the qualifier \emph{stable} if
$\alpha_i,\alpha_j\geq 0$.
When $\mathbf q$ is clear from context, we omit
``relative to $\mathbf q$.''
\end{definition}

In our applications, $\mathbf q$ is specified by the
relevant neighborhood decomposition.
The coefficients $\alpha_i,\alpha_j$ are effective
coefficients of these relations and need not coincide
with the original node strengths $\mu_i,\mu_j$.

\subsection{Structural twins at equilibrium}
Structural twins provide the most direct instances of the geometric twin relations introduced above.
\begin{corollary}\label{closed-twins}
Let $\boldsymbol{\theta}$ be a first-order stationary point of $E_G$ in~\eqref{E_2}. Then structural closed twins $a$ and $b$ must have phasors satisfying one of the three cases for \eqref{acase1} listed in Lemma~\ref{geo}, with $\boldsymbol q=\sum_{j \in N(a)\setminus\{b\}} \mathbf{v}_j = \sum_{j \in N(b)\setminus\{a\}} \mathbf{v}_j$.
\end{corollary}

The cases in Lemma~5.1 are illustrated in Figure~\ref{fig:closed-twins}, 
where the synchronized case is split into two parameter regimes.

\begin{proof}
Since $\boldsymbol{\theta}$ is a first-order stationary point of \eqref{E_2}, according to Lemma \ref{eq-vec}, there exist $\mu_a,\mu_b \in \mathbb R$ such that $\mu_a\mathbf{v}_a = \sum_{j\in N(a)} \mathbf{v}_j$, and $\mu_b\mathbf{v}_b=\sum_{j\in N(b)} \mathbf{v}_j$.
Let $\mathbf{q}=\sum_{j \in N(a)\setminus\{b\}} \mathbf{v}_j$, then equivalently $\mu_a\mathbf{v}_a = \mathbf{q} + \mathbf{v}_b$, and $\mu_b\mathbf{v}_b = \mathbf{q} + \mathbf{v}_a$.
Thus the conclusion holds according to Lemma \ref{geo}.
\end{proof}

\begin{corollary}\label{stable-closed-twins}
Let $\boldsymbol{\theta}$ be a second-order stationary point of $E_G$ in~\eqref{E_2}.
Then structural closed twins $a$ and $b$ must synchronize, i.e., $\theta_a=\theta_b \pmod{2\pi}$.
\end{corollary}

\begin{proof}
   Lemma \ref{stable-eq} tells both $\mu_a$ and $\mu_b$ are nonnegative, since $\boldsymbol{\theta}$ is a second-order stationary point of \eqref{E_2}. According to Corollary \ref{closed-twins}, $\mathbf{v}_a=\mathbf{v}_b$.
\end{proof}

\begin{figure}[H]
\centering
\begin{tikzpicture}[
    x=1cm, y=1cm, font=\small,
    ct edge/.style={
        draw=black!60,
        line width=0.45pt
    },
    ct twin/.style={
        circle, draw=black, fill=black,
        minimum size=1.6mm,
        inner sep=0pt, outer sep=0pt
    },
    ct neighbor/.style={
        circle, draw=black, fill=white,
        line width=0.4pt,
        minimum size=1.6mm,
        inner sep=0pt, outer sep=0pt
    },
    ct stub/.style={
        draw=black!28,
        line width=0.35pt
    },
    ct vector/.style={
        line width=0.75pt,
        -{Stealth[length=1.8mm,width=1.25mm]}
    },
    ct sum/.style={
        line width=1pt,
        -{Stealth[length=2mm,width=1.4mm]}
    },
    ct closing/.style={
        draw=red!80!black,
        densely dashed,
        line width=0.65pt,
        -{Stealth[length=1.7mm,width=1.15mm]}
    },
    ct dot/.style={
        circle, fill=black,
        inner sep=0pt,
        minimum size=1.25mm
    },
    ct positive/.style={
        draw=green!55!black,
        fill=green!55!black
    },
    ct negative/.style={
        draw=red!80!black,
        fill=red!80!black
    },
    every node/.style={inner sep=1.5pt}
]

\foreach \xx/\acolor/\bcolor in {
    0/ct positive/ct positive,
    4/ct negative/ct negative,
    8/ct negative/ct positive,
    12/ct negative/ct negative%
} {
    \begin{scope}[xshift=\xx cm]
        \coordinate (a) at (-0.34,2.58);
        \coordinate (b) at ( 0.34,2.58);

        \draw[ct edge] (a) -- (b);

        \foreach \k/\px in {
            1/-1.00,
            2/-0.60,
            3/-0.20,
            4/0.20,
            5/0.60,
            6/1.00%
        } {
            \coordinate (n\k) at (\px,1.86);
            \draw[ct edge] (a) -- (n\k) -- (b);

            \foreach \dx in {-0.14,0,0.14} {
                \draw[ct stub]
                    (n\k) -- ++(\dx,-0.24);
            }
        }

        \foreach \k in {1,2,3,4,5,6} {
            \node[ct neighbor] at (n\k) {};
        }

        \node[ct twin,\acolor] at (a) {};
        \node[ct twin,\bcolor] at (b) {};

        \draw[line width=0.55pt]
            (0,0) circle[radius=0.92];
        \fill (0,0) circle[radius=0.65pt];
    \end{scope}
}

\begin{scope}
    \draw[ct sum] (0,0) -- (1.53,0);
    \node[ct dot,ct positive] at (0.92,0) {};
\end{scope}

\begin{scope}[xshift=4cm]
    \draw[ct sum] (0,0) -- (1.53,0);
    \node[ct dot,ct negative] at (-0.92,0) {};
\end{scope}

\begin{scope}[xshift=8cm]
    \draw[ct sum] (0,0) -- (1.53,0);
    \node[ct dot,ct negative] at (-0.92,0) {};
    \node[ct dot,ct positive] at ( 0.92,0) {};
\end{scope}

\begin{scope}[xshift=12cm]
    \coordinate (O) at (0,0);
    \coordinate (L) at (-0.92,0);
    \coordinate (U) at (-0.736,0.552);
    \coordinate (R) at (0.92,0);

    \draw[dashed,line width=0.3pt] (O) -- (R);
    \draw[dashed,ct vector] (O) -- (L);
    \draw[dashed,ct vector] (O) -- (U);
    \draw[ct closing,draw=black,solid] (U) -- (R);

    \node[font=\scriptsize] at (0.02,0.57)
        {$\boldsymbol{q}$};

    \node[ct dot,ct negative] at (L) {};
    \node[ct dot,ct negative] at (U) {};
    \node[
        ct dot,
        draw=black,
        fill=white,
        line width=0.4pt
    ] at (R) {};
\end{scope}

\node at (0,-1.45)
    {$\mu_a=\mu_b>-1$};

\node at (4,-1.45)
    {$\mu_a=\mu_b<-1$};

\node at (8,-1.45)
    {$\mu_a+\mu_b=-2,\ \mu_a\ne-1$};

\node at (12,-1.45)
    {$\mu_a=\mu_b=-1$};

\foreach \xx/\panelletter in {0/a,4/b,8/c,12/d} {
    \node at (\xx,-2.10) {(\panelletter)};
}

\end{tikzpicture}

\caption{Representative phasor configurations of structural
closed twins at equilibrium. Green indicates $\mu>0$,
and red indicates $\mu<0$. In panels (a)--(c), the arrow indicates
the direction of $\mathbf q$.}
\label{fig:closed-twins}
\end{figure}

\begin{corollary}\label{coro:structural-open-twins}
Let $\boldsymbol{\theta}$ be a first-order stationary point of $E_G$ in~\eqref{E_2}.
If $a$ and $b$ are structural open twins with
\[
N(a)=N(b)=N_{\mathrm{common}},
\]
then their phasors satisfy one of the three cases for~\eqref{acase2} listed in Lemma~\ref{geo}, with
\[
\boldsymbol q=\sum_{j\in N_{\mathrm{common}}}\boldsymbol v_j.
\]
\end{corollary}

\begin{proof}
Since $\boldsymbol{\theta}$ is a first-order stationary point of \eqref{E_2}, according to Lemma \ref{eq-vec}, there exist $\mu_a,\mu_b \in \mathbb R$ such that
\begin{equation}
    \mu_a\mathbf{v}_a = \sum_{j\in N(a)} \mathbf{v}_j, \quad \mathrm{and} \quad \mu_b\mathbf{v}_b=\sum_{j\in N(b)} \mathbf{v}_j.
\end{equation}
Since $a$ and $b$ form a pair of structural open twins, let $\boldsymbol q=\sum_{j \in N(a)} \mathbf{v}_j=\sum_{j \in N(b)} \mathbf{v}_j$, then equivalently $\mu_a\mathbf{v}_a = \boldsymbol q$ and $\mu_b\mathbf{v}_b = \boldsymbol q$.
Then the conclusion holds according to Lemma \ref{geo}.
\end{proof}

\begin{corollary}\label{stable-open-twins}
Let $\boldsymbol{\theta}$ be a second-order stationary point of the Kuramoto model on $G$.
If $a$ and $b$ are structural open twins with
\[
N(a)=N(b)=N_{\mathrm{common}},
\]
then they synchronize whenever
\[
\sum_{i\in N_{\mathrm{common}}}\boldsymbol v_i\neq \boldsymbol 0.
\]
\end{corollary}

\begin{proof}
   Lemma \ref{stable-eq} tells both $\mu_a$ and $\mu_b$ are nonnegative, since $\boldsymbol{\theta}$ is a second-order stationary point of \eqref{E_2}. Since $\sum_{i \in N_{\mathrm{common}}} \mathbf{v}_i \neq \mathbf{0}$, then $\mathbf{v}_a=\mathbf{v}_b$ according to Corollary \ref{coro:structural-open-twins}.
\end{proof}

\begin{figure}[H]
\centering
\begin{tikzpicture}[
    x=1cm, y=1cm, font=\small,
    ot edge/.style={draw=black!60,line width=0.45pt},
    ot twin/.style={
        circle,draw=black,fill=black,
        minimum size=1.6mm,
        inner sep=0pt,outer sep=0pt
    },
    ot neighbor/.style={
        circle,draw=black,fill=white,
        line width=0.4pt,
        minimum size=1.6mm,
        inner sep=0pt,outer sep=0pt
    },
    ot stub/.style={draw=black!28,line width=0.35pt},
    ot sum/.style={
        line width=1pt,
        -{Stealth[length=2mm,width=1.4mm]}
    },
    ot dot/.style={
        circle,fill=black,
        inner sep=0pt,minimum size=1.25mm
    },
    ot positive/.style={
        draw=green!55!black,fill=green!55!black
    },
    ot negative/.style={
        draw=red!80!black,fill=red!80!black
    },
    ot zero/.style={draw=black,fill=black},
    every node/.style={inner sep=1.5pt}
]

\foreach \xx/\acolor/\bcolor in {
    0/ot positive/ot positive,
    4/ot negative/ot negative,
    8/ot negative/ot positive,
    12/ot zero/ot zero%
} {
    \begin{scope}[xshift=\xx cm]
        \coordinate (a) at (-0.34,2.58);
        \coordinate (b) at ( 0.34,2.58);

        \foreach \k/\px in {
            1/-1.00,2/-0.60,3/-0.20,
            4/0.20,5/0.60,6/1.00%
        } {
            \coordinate (n\k) at (\px,1.86);
            \draw[ot edge] (a) -- (n\k) -- (b);

            \foreach \dx in {-0.14,0,0.14} {
                \draw[ot stub]
                    (n\k) -- ++(\dx,-0.24);
            }
        }

        \foreach \k in {1,2,3,4,5,6} {
            \node[ot neighbor] at (n\k) {};
        }

        \node[ot twin,\acolor] at (a) {};
        \node[ot twin,\bcolor] at (b) {};

        \draw[line width=0.55pt]
            (0,0) circle[radius=0.92];
    \end{scope}
}

\begin{scope}
    \draw[ot sum] (0,0) -- (1.53,0);
    \node[ot dot,ot positive] at (0.92,0) {};
\end{scope}

\begin{scope}[xshift=4cm]
    \draw[ot sum] (0,0) -- (1.53,0);
    \node[ot dot,ot negative] at (-0.92,0) {};
\end{scope}

\begin{scope}[xshift=8cm]
    \draw[ot sum] (0,0) -- (1.53,0);
    \node[ot dot,ot negative] at (-0.92,0) {};
    \node[ot dot,ot positive] at ( 0.92,0) {};
\end{scope}

\begin{scope}[xshift=12cm]
    \node[ot dot,ot zero] at (0.92,0) {};
    \node[ot dot,ot zero]
        at ({0.92*cos(35)},{-0.92*sin(35)}) {};
\end{scope}

\node at (0,-1.45)
    {$\mu_a=\mu_b>0$};

\node at (4,-1.45)
    {$\mu_a=\mu_b<0$};

\node at (8,-1.45)
    {$\mu_a=-\mu_b\ne0$};

\node at (12,-1.45)
    {$\mu_a=\mu_b=0$};

\foreach \xx/\panelletter in {0/a,4/b,8/c,12/d} {
    \node at (\xx,-2.10) {(\panelletter)};
}

\end{tikzpicture}

\caption{Representative phasor configurations of structural
open twins at equilibrium. Green indicates $\mu>0$,
red indicates $\mu<0$, and black indicates $\mu=0$.}
\label{fig:open-twins}
\end{figure}

As a first application, we obtain the following class of second-order globally synchronizing graphs.

\begin{corollary}
Complete split graphs are second-order globally synchronizing. 
\end{corollary}

\begin{proof}
Let $\boldsymbol{\theta}$ be a second-order stationary point of $E_G$ on a complete split graph $G$ with vertex partition
\[
V=C\sqcup I,
\]
where $C$ is a nonempty clique, $I$ is an independent set, and every node in $I$ is adjacent to every node in $C$.

We first show that all nodes in $C$ synchronize. Indeed, for any $u,v\in C$, we have
\[
N[u]=N[v]=V,
\]
so $u$ and $v$ are structural closed twins. Hence, by Corollary~\ref{stable-closed-twins}, they synchronize at $\boldsymbol{\theta}$.
Let $\mathbf{v}_C$ denote the common phasor of nodes in $C$. Then
\[
\sum_{j\in C}\mathbf{v}_j
=
|C|\,\mathbf{v}_C
\neq \mathbf{0}.
\]
If $I=\varnothing$, then $V=C$, so all nodes of $G$
are already synchronized.
Hence assume that $I\neq\varnothing$.
Next, for any $u,v\in I$, we have
\[
N(u)=N(v)=C,
\]
so $u$ and $v$ are structural open twins with common neighborhood $C$. Since
\[
\sum_{j\in C}\mathbf{v}_j
=
|C|\,\mathbf{v}_C
\neq \mathbf{0},
\]
Corollary~\ref{stable-open-twins} implies that all nodes in $I$ synchronize. Let $\mathbf{v}_I$ denote their common phasor.

Finally, fix any node $i\in I$. Since $\boldsymbol{\theta}$ is a second-order stationary point, Lemma~\ref{stable-eq} implies that
\[
\sum_{j\in N(i)}\mathbf{v}_j
=
\mu_i\mathbf{v}_i,
\qquad
\mu_i\ge 0.
\]
Since $N(i)=C$, we obtain
\[
|C|\,\mathbf{v}_C
=
\sum_{j\in C}\mathbf{v}_j
=
\mu_i\mathbf{v}_i.
\]
The left-hand side is nonzero, so $\mu_i>0$. Since both $\mathbf{v}_C$ and $\mathbf{v}_i$ are unit vectors, it follows that
\[
\mathbf{v}_i=\mathbf{v}_C.
\]
As $i\in I$ was arbitrary, all nodes in $I$ synchronize with the clique $C$. Hence all nodes of $G$ have the same phasor, so $\boldsymbol{\theta}$ is fully synchronized. Therefore, every complete split graph is second-order globally synchronizing.
\end{proof}
\begin{figure}[H]
\centering
\begin{tikzpicture}[scale=1, line cap=round, line join=round]

\def\rin{1.0}
\def\rout{3.0}

\foreach[count=\k from 1] \ang in {90,30,-30,-90,-150,150}{
    \coordinate (m\k) at (\ang:\rin);
}

\foreach[count=\k from 1] \ang in {90,45,0,-45,-90,-135,180,135}{
    \coordinate (p\k) at (\ang:\rout);
}

\foreach \i in {1,...,8}{
    \foreach \j in {1,...,6}{
        \draw[black, line width=0.5pt] (p\i)--(m\j);
    }
}

\foreach \u/\v in {
    m1/m2, m1/m3, m1/m4, m1/m5, m1/m6,
    m2/m3, m2/m4, m2/m5, m2/m6,
    m3/m4, m3/m5, m3/m6,
    m4/m5, m4/m6,
    m5/m6}{
    \draw[black, line width=0.5pt] (\u)--(\v);
}

\foreach \p in {m1,m2,m3,m4,m5,m6,p1,p2,p3,p4,p5,p6,p7,p8}{
    \fill (\p) circle (2.2pt);
}
\end{tikzpicture}
\caption{A complete split graph.}
\end{figure}
\subsection{Benign extra neighbors}\label{sec:benign-extra}
This section presents a setting in which geometrically stable closed twins arise even when the corresponding nodes are not structural closed twins. See Figure~\ref{fig:lemma5_3}
for an illustration.

\begin{lemma}\label{pairwise-no-nauty-extra-node}
Let $\boldsymbol{\theta}$ be a second-order stationary point of $E_G$ on a graph
$G=(V,E)$.
Let $a,b\in V$ be adjacent nodes, and suppose that there exist disjoint sets
$S,T\subseteq V\setminus\{a,b\}$ such that
\[
N(a)\setminus\{b\}=S\sqcup T,
\qquad
N(b)\setminus\{a\}=S.
\]
If
\[
\theta_i=\theta_b \pmod{2\pi}, \qquad\text{for every } i\in T,
\]
then $a$ and $b$ synchronize, i.e.,
\[
\theta_a=\theta_b \pmod{2\pi}.
\]
\end{lemma}

\begin{proof}
Since $\boldsymbol{\theta}$ is a second-order stationary point,
Lemma~\ref{stable-eq} implies that there exist $\mu_a,\mu_b\ge 0$ such that
\[
\mu_a\mathbf{v}_a
=
\mathbf{v}_b+\sum_{j\in S\sqcup T}\mathbf{v}_j,
\qquad
\mu_b\mathbf{v}_b
=
\mathbf{v}_a+\sum_{j\in S}\mathbf{v}_j.
\]
Let
\[
\mathbf q=\sum_{j\in S\sqcup T}\mathbf{v}_j.
\]
Since $\mathbf{v}_i=\mathbf{v}_b$ for every $i\in T$, we have
\[
\mathbf q
=
\sum_{j\in S}\mathbf{v}_j+|T|\mathbf{v}_b.
\]
Hence the two equilibrium equations can be rewritten as
\[
\mu_a\mathbf{v}_a
=
\mathbf{v}_b+\mathbf q,
\qquad
(\mu_b+|T|)\mathbf{v}_b
=
\mathbf{v}_a+\mathbf q.
\]
Thus $a$ and $b$ form a pair of geometrically closed twins with effective coefficients
$\mu_a$ and $\mu_b+|T|$. Since
\[
\mu_a\ge0
\qquad\text{and}\qquad
\mu_b+|T|\ge0,
\]
the second and third cases in Lemma~\ref{geo} are impossible. Hence only the
synchronized case remains, and therefore
\[
\mathbf{v}_a=\mathbf{v}_b.
\]
Thus $a$ and $b$ synchronize, i.e., \[
\theta_a=\theta_b \pmod{2\pi}.
\]
\end{proof}

\begin{figure}[H]
    \centering
\begin{tikzpicture}[
    x=1cm,y=1cm,font=\small,
    vertex/.style={
        circle,draw=black,fill=white,
        line width=0.45pt,
        minimum size=1.6mm,
        inner sep=0pt,outer sep=0pt
    },
    highlighted/.style={
        vertex,fill=green!55!black
    },
    edge/.style={line width=0.7pt},
    region/.style={fill=black!16,draw=none},
    every label/.style={inner sep=1pt},
    label distance=2pt
]

\foreach \xx/\stage in {0/0,6.8/1} {
    \begin{scope}[xshift=\xx cm]

        \ifnum\stage=0\relax
            \path[region]
                (-0.12,-0.19)
                .. controls (-0.36,-0.20) and (-0.37,0.10)
                    .. (-0.08,0.28)
                .. controls (0.40,0.55) and (1.30,0.86)
                    .. (1.61,1.37)
                .. controls (1.66,1.70) and (2.10,1.70)
                    .. (2.08,1.43)
                .. controls (2.07,1.18) and (1.36,0.91)
                    .. (0.91,0.62)
                .. controls (0.72,0.43) and (0.86,0.10)
                    .. (0.78,-0.10)
                .. controls (0.70,-0.23) and (0.08,-0.20)
                    .. cycle;
        \else
            \path[region]
                (-0.12,-0.19)
                .. controls (-0.37,-0.20) and (-0.37,0.08)
                    .. (-0.19,0.33)
                .. controls (0.00,0.64) and (0.24,0.98)
                    .. (0.36,1.42)
                .. controls (0.41,1.66) and (0.52,1.73)
                    .. (0.78,1.68)
                .. controls (1.06,1.61) and (1.55,1.65)
                    .. (1.84,1.67)
                .. controls (2.13,1.71) and (2.15,1.35)
                    .. (1.94,1.25)
                .. controls (1.61,1.09) and (1.05,0.93)
                    .. (0.85,0.60)
                .. controls (0.67,0.35) and (0.87,0.12)
                    .. (0.79,-0.09)
                .. controls (0.72,-0.23) and (0.02,-0.22)
                    .. cycle;
        \fi

        \coordinate (A) at (0.65,1.45);
        \coordinate (B) at (1.85,1.45);
        \coordinate (T1) at (-0.10,0);
        \coordinate (T2) at (0.65,0);
        \coordinate (S) at (1.27,0);

        \foreach \u/\v in {A/B,A/T1,A/T2,A/S,B/S} {
            \draw[edge] (\u) -- (\v);
        }

        \node[highlighted,label=above:{$a$}] at (A) {};
        \node[highlighted,label=above:{$b$}] at (B) {};
        \node[vertex] at (T1) {};
        \node[vertex] at (T2) {};
        \node[vertex] at (S) {};

        \node[inner sep=0pt] at (0.275,0) {$\cdots$};

        \draw[
            decorate,
            decoration={brace,mirror,amplitude=3pt},
            line width=0.45pt
        ] (-0.23,-0.29) -- (0.78,-0.29);

        \node at (0.275,-0.62) {$T$};
        \node at (1.27,-0.62) {$S$};

    \end{scope}
}

\draw[
    line width=0.6pt,
    -{Stealth[length=2.4mm,width=1.8mm]}
] (3.15,0.66) -- (5.15,0.66);

\end{tikzpicture}

\caption{Illustration of Lemma~\ref{pairwise-no-nauty-extra-node}.
Nodes $a$ and $b$ share the common neighbors $S$, while $a$
has additional neighbors $T$. Gray regions indicate synchronized
nodes. Under the assumptions of the lemma, synchronization of
$b$ with all nodes in $T$ implies that $a$ also synchronizes
with them.}
    \label{fig:lemma5_3}
\end{figure}

\subsection{Synchronous homogeneous extension}
This section considers synchronous homogeneous extensions,
where synchronized nodes share a common external neighborhood
but may have different internal connections.
At any second-order stationary point, the sum of their external
neighbors' vectors is a nonnegative multiple of their common vector. See Figure~\ref{fig:synchronous-extension}
for an illustration.
\begin{lemma}\label{lem:opentwins}
Let $\boldsymbol\theta$ be a second-order stationary point of the energy
function~(2) of the Kuramoto model on a graph $G=(V,E)$.
Let $P,Q\subseteq V$ be disjoint nonempty sets,
and suppose that
\[
    P\subseteq N(i)\subseteq P\cup Q
    \qquad \text{for every } i\in Q.
\]
If the nodes in $Q$ are synchronized, with
$\mathbf{v}_i=\mathbf{v}$ for all $i\in Q$, then there exists
a scalar $\mu\geq 0$ such that
\[
    \sum_{j\in P}\mathbf{v}_j=\mu\mathbf{v}.
\]
\end{lemma}

\begin{proof}[Proof of Lemma~\ref{lem:opentwins}]
Fix $i_0\in Q$ and let $Q_{i_0}=N(i_0)\cap Q$.
By assumption, $N(i_0)=P\sqcup Q_{i_0}$.
First-order stationarity implies that
\[
    \sum_{j\in N(i_0)}\mathbf{v}_j
    =\lambda\mathbf{v}_{i_0}
\]
for some $\lambda\in\mathbb{R}$.
Since $\mathbf{v}_i=\mathbf{v}$ for all $i\in Q$, it follows that
\[
    \sum_{j\in P}\mathbf{v}_j
    =(\lambda-|Q_{i_0}|)\mathbf{v}
    =\mu\mathbf{v},
\]
where $\mu:=\lambda-|Q_{i_0}|$.

It remains to show that $\mu\ge0$.
Let $H$ be the Hessian of the energy function~\eqref{E_2}
at $\boldsymbol{\theta}$, and define
$\boldsymbol{x}\in\mathbb{R}^{|V|}$ by
\[
    x_k=
    \begin{cases}
        1, & k\in Q,\\
        0, & k\notin Q.
    \end{cases}
\]
Since $P\subseteq N(i)\subseteq P\cup Q$ for every $i\in Q$,
the edges crossing from $Q$ to $V\setminus Q$ are precisely
those joining $Q$ to $P$.
Thus, by second-order stationarity,
\[
    0\le \boldsymbol{x}^{\top}H\boldsymbol{x}
    =\sum_{i\in Q}\sum_{j\in P}\cos(\theta_i-\theta_j)
    =|Q|\left\langle
        \mathbf{v},\sum_{j\in P}\mathbf{v}_j
      \right\rangle
    =|Q|\mu,
\]
where the last equality uses $\|\mathbf{v}\|=1$.
Since $Q\neq\varnothing$, we conclude that $\mu\ge0$.
\end{proof}

\begin{figure}[H]
\centering
\begin{tikzpicture}[
    x=1cm,y=1cm,font=\small,
    se vertex/.style={circle,draw=black,fill=green!55!black,
        line width=0.4pt,minimum size=1.6mm,
        inner sep=0pt,outer sep=0pt},
    se set/.style={circle,draw=black,fill=white,
        line width=0.45pt,minimum size=4mm,
        inner sep=0pt,outer sep=0pt},
    se edge/.style={draw=black,line width=0.45pt},
    se internal/.style={draw=black!55,line width=0.35pt},
    se emphasis/.style={draw=black!24,line width=2.8pt,line cap=round},
    se region/.style={fill=black!16,draw=none}
]
    \foreach \xx/\stage in {0/0,4.1/1} {
        \begin{scope}[xshift=\xx cm]
            \coordinate (P) at (0,2.10);
            \coordinate (Q1) at (-0.40,0.48);
            \coordinate (Q2) at (-0.98,-0.05);
            \coordinate (Q3) at (0,-0.65);
            \coordinate (Q4) at (0.40,0.48);
            \coordinate (Q5) at (0.98,-0.05);

            \path[se region]
                (0,-0.06) ellipse[x radius=1.20,y radius=0.82];

            \foreach \u/\v in {1/2,2/3,3/5,5/4,4/1,1/5,2/4} {
                \draw[se internal] (Q\u) -- (Q\v);
            }

            \ifnum\stage=1\relax
                \foreach \k in {1,2,3,4,5} {
                    \draw[se emphasis] (P) -- (Q\k);
                }
            \fi

            \foreach \k in {1,2,3,4,5} {
                \draw[se edge] (P) -- (Q\k);
            }

            \node[se set] at (P) {};
            \foreach \k in {1,2,3,4,5} {
                \node[se vertex] at (Q\k) {};
            }
        \end{scope}
    }

    \node at (-0.65,2.10) {$P$};
    \node at (-1.45,-0.05) {$Q$};

    \draw[line width=0.6pt,
        -{Stealth[length=2mm,width=1.4mm]}]
        (1.45,1.10) -- (2.65,1.10);
\end{tikzpicture}
\caption{Illustration of Lemma~\ref{lem:opentwins}.
The shaded set $Q$ is synchronized and has common external
neighborhood $P$, represented by the large circle.
At a second-order stationary point,
$\sum_{j\in P}\mathbf{v}_j=\mu\mathbf{v}$ for some $\mu\ge0$,
where $\mathbf{v}$ is the common phasor of the nodes in $Q$.}
\label{fig:synchronous-extension}
\end{figure}

\section{Propagation argument on rooted tree:  proof of Theorem \ref{main1}}\label{sec:quasi-threshold}
\subsection{Leaf-like node and its upward propagation}
We introduce leaf-like nodes and show how the leaf-like property
propagates upward through the underlying rooted tree.
\begin{definition}[Leaf-like node in $T$ at state $\boldsymbol{\theta}$]
Let $G=(V,E)$ be a connected quasi-threshold graph with
underlying rooted tree $T$, and let
$\boldsymbol{\theta}\in\mathbb{R}^{|V|}$ be a state of the
Kuramoto model on $G$.
For $l\in V(T)$, let $\mathrm{Desc}_T(l)$ denote the set of
descendants of $l$ in $T$.
We say that $l$ is \emph{leaf-like in $T$ at
$\boldsymbol{\theta}$} if
\[
    \theta_i \equiv \theta_l \pmod{2\pi}
    \qquad \text{for every } i\in\mathrm{Desc}_T(l).
\]
\end{definition}

\begin{lemma}[Propagation of the leaf-like property]
\label{lemma:leaflike}
Let $\boldsymbol{\theta}\in\mathbb{R}^{|V|}$ be a second-order
stationary point of the Kuramoto energy $E_G(\boldsymbol{\theta})$ on a connected
quasi-threshold graph $G=(V,E)$ with underlying rooted tree $T$.
If $a\in V(T)$ is a non-leaf node and all its children are
leaf-like in $T$ at $\boldsymbol{\theta}$, then $a$ is also
leaf-like in $T$ at $\boldsymbol{\theta}$.
\end{lemma}

\begin{proof}[Proof of Lemma~\ref{lemma:leaflike}]
If $a$ has exactly one child $b$, then
\[
    N_G[a]=N_G[b]
    =\mathrm{Anc}_T(a)\cup\{a,b\}\cup\mathrm{Desc}_T(b).
\]
Thus, $a$ and $b$ are structural closed twins in $G$ and
synchronize at $\boldsymbol{\theta}$ by
Corollary~\ref{stable-closed-twins}.
Since $b$ is leaf-like in $T$ at $\boldsymbol{\theta}$ and
\[
    \mathrm{Desc}_T(a)=\{b\}\cup\mathrm{Desc}_T(b),
\]
every descendant of $a$ synchronizes with $a$.
Hence, $a$ is leaf-like in $T$ at $\boldsymbol{\theta}$.
It remains to consider the case where $a$ has at least two
children.

\medskip
\noindent
\textit{Step 1: Applying the synchronous extension lemma.}
Set
\[
    P:=\{a\}\cup\mathrm{Anc}_T(a),
    \qquad
    \mathbf{s}:=\sum_{j\in P}\mathbf{v}_j
    =\mathbf{v}_a+\sum_{j\in\mathrm{Anc}_T(a)}\mathbf{v}_j.
\]
Fix a child $i$ of $a$ and let
\[
    Q_i:=\{i\}\cup\mathrm{Desc}_T(i).
\]
Since $i$ is leaf-like in $T$ at $\boldsymbol{\theta}$,
we have $\mathbf{v}_k=\mathbf{v}_i$ for every $k\in Q_i$.
Moreover, adjacency in $G$ is determined by the
ancestor--descendant relation in $T$.
Thus, every node $k\in Q_i$ is adjacent to every node in $P$
and has no neighbors outside $P\cup Q_i$; that is,
\[
    P\subseteq N_G(k)\subseteq P\cup Q_i
    \qquad \text{for every } k\in Q_i.
\]
The sets $P$ and $Q_i$ are disjoint and nonempty.
Applying Lemma~\ref{lem:opentwins} therefore gives a scalar
$\alpha_i\ge0$ such that
\[
    \mathbf{s}=\alpha_i\mathbf{v}_i.
\]
Consequently, if $\mathbf{s}\neq\boldsymbol{0}$, then
$\alpha_i=\|\mathbf{s}\|>0$ and
\[
    \mathbf{v}_i=\frac{\mathbf{s}}{\|\mathbf{s}\|}
    \qquad \text{for every } i\in\mathrm{Children}_T(a).
\]
Hence all children of $a$ synchronize whenever
$\mathbf{s}\neq\boldsymbol{0}$.

\medskip
\noindent
\textit{Step 2: Excluding a vanishing external sum.}
Suppose, for contradiction, that $\mathbf{s}=\boldsymbol{0}$.
Write
\[
    \mathbf{u}:=\sum_{j\in\mathrm{Anc}_T(a)}\mathbf{v}_j,
    \qquad
    \mathbf{d}:=\sum_{j\in\mathrm{Desc}_T(a)}\mathbf{v}_j.
\]
Then $\mathbf{u}=-\mathbf{v}_a$.
Since
\[
    N_G(a)=\mathrm{Anc}_T(a)\cup\mathrm{Desc}_T(a),
\]
Lemma~\ref{stable-eq} gives
\[
    \mathbf{u}+\mathbf{d}=\lambda\mathbf{v}_a
    \qquad \text{for some } \lambda\ge0.
\]
It follows that
\[
    \mathbf{d}=(\lambda+1)\mathbf{v}_a.
\]

Let $H=\nabla^2 E_G(\boldsymbol{\theta})$ and define
$\boldsymbol{x}\in\mathbb{R}^{|V|}$ by
\[
    x_k=
    \begin{cases}
        1, & k\in\{a\}\cup\mathrm{Desc}_T(a),\\
        0, & \text{otherwise}.
    \end{cases}
\]
By the ancestor--descendant structure of $G$, the edges
crossing from $\{a\}\cup\mathrm{Desc}_T(a)$ to its complement
are precisely all edges joining this set to
$\mathrm{Anc}_T(a)$.
Therefore,
\[
\begin{aligned}
    \boldsymbol{x}^{\top}H\boldsymbol{x}
    &=
    \sum_{k\in\{a\}\cup\mathrm{Desc}_T(a)}
    \sum_{j\in\mathrm{Anc}_T(a)}
    \cos(\theta_k-\theta_j)\\
    &=\langle\mathbf{v}_a+\mathbf{d},\mathbf{u}\rangle\\
    &=\langle(\lambda+2)\mathbf{v}_a,-\mathbf{v}_a\rangle\\
    &=-(\lambda+2)<0,
\end{aligned}
\]
where we used $\|\mathbf{v}_a\|=1$.
This contradicts the positive semidefiniteness of $H$.
Thus $\mathbf{s}\neq\boldsymbol{0}$, and Step~1 implies
that all children of $a$ synchronize.
Since each child is leaf-like in $T$ at
$\boldsymbol{\theta}$, all descendants of $a$ synchronize
with one another.

\medskip
\noindent
\textit{Step 3: Synchronizing $a$ with its descendants.}
Fix a child $b$ of $a$ and set
\[
\begin{aligned}
    S&:=\mathrm{Anc}_T(a)\cup\mathrm{Desc}_T(b),\\
    R&:=\mathrm{Desc}_T(a)\setminus
    \bigl(\{b\}\cup\mathrm{Desc}_T(b)\bigr).
\end{aligned}
\]
These are disjoint subsets of $V\setminus\{a,b\}$ satisfying
\[
    N_G(a)\setminus\{b\}=S\sqcup R,
    \qquad
    N_G(b)\setminus\{a\}=S.
\]
Every node in $R$ synchronizes with $b$, because all
descendants of $a$ are synchronized.
Applying Lemma~\ref{pairwise-no-nauty-extra-node} to the
adjacent nodes $a,b$, with common neighbors $S$ and extra
neighbors $R$, yields $\mathbf{v}_a=\mathbf{v}_b$.
Thus every descendant of $a$ synchronizes with $a$,
so $a$ is leaf-like in $T$ at $\boldsymbol{\theta}$.
\end{proof}

\begin{figure}[H]
\centering
\begin{tikzpicture}[
    x=1cm,y=1cm,font=\small,
    up vertex/.style={circle,draw=black,fill=green!55!black,
        line width=0.4pt,minimum size=1.6mm,
        inner sep=0pt,outer sep=0pt},
    up set/.style={circle,draw=black,fill=white,
        line width=0.45pt,minimum size=4mm,
        inner sep=0pt,outer sep=0pt},
    up tree/.style={draw=black,line width=0.55pt},
    up ancestor/.style={draw=black!38,line width=0.35pt},
    up region/.style={fill=black!13,draw=none},
    up edge mark/.style={draw=orange,line width=0.90pt,
        line cap=round},
    every node/.style={inner sep=1pt}
]
    \foreach \xx/\yy/\stage in {
        0/0/1,5.4/0/2,0/-4.75/3,5.4/-4.75/4%
    } {
        \begin{scope}[shift={(\xx,\yy)}]
            \coordinate (P) at (0,2.80);
            \coordinate (A) at (0,1.70);
            \coordinate (B) at (-1.15,0.80);
            \coordinate (C) at (-0.38,0.80);
            \coordinate (D) at (0.38,0.80);
            \coordinate (E) at (1.15,0.80);
            \coordinate (BL) at (-1.70,0.15);
            \coordinate (BR) at (-0.60,0.15);
            \coordinate (EL) at (0.60,0.15);
            \coordinate (ER) at (1.70,0.15);

            \ifnum\stage<3\relax
                \foreach \cx in {-1.15,1.15} {
                    \path[up region]
                        ({\cx-0.87},-0.43)
                        .. controls ({\cx-1.00},-0.18)
                            and ({\cx-0.45},0.59)
                            .. ({\cx-0.12},0.98)
                        .. controls ({\cx-0.04},1.10)
                            and ({\cx+0.08},1.10)
                            .. ({\cx+0.19},0.92)
                        .. controls ({\cx+0.49},0.49)
                            and ({\cx+0.96},-0.16)
                            .. ({\cx+0.88},-0.38)
                        .. controls ({\cx+0.80},-0.58)
                            and ({\cx-0.57},-0.56)
                            .. cycle;
                }
            \else
                \ifnum\stage=3\relax
                    \path[up region]
                        (-2.00,-0.43)
                        .. controls (-2.18,-0.45) and (-1.90,0.35)
                            .. (-1.45,0.90)
                        .. controls (-1.10,1.20) and (-0.40,1.25)
                            .. (0,1.20)
                        .. controls (0.40,1.25) and (1.10,1.20)
                            .. (1.45,0.90)
                        .. controls (1.90,0.35) and (2.18,-0.45)
                            .. (2.00,-0.43)
                        .. controls (1.00,-0.60) and (-1.00,-0.60)
                            .. cycle;
                \else
                    \path[up region]
                        (-2.00,-0.43)
                        .. controls (-2.15,-0.36) and (-1.87,0.24)
                            .. (-1.30,1.04)
                        .. controls (-0.92,1.55) and (-0.35,2.05)
                            .. (0,2.02)
                        .. controls (0.35,2.05) and (0.92,1.55)
                            .. (1.30,1.04)
                        .. controls (1.87,0.24) and (2.15,-0.36)
                            .. (2.00,-0.43)
                        .. controls (1.00,-0.67) and (-1.00,-0.67)
                            .. cycle;
                \fi
            \fi

            \foreach \target/\tx/\ty in {
                A/0/1.70,
                B/-1.15/0.80,C/-0.38/0.80,
                D/0.38/0.80,E/1.15/0.80,
                BL/-1.70/0.15,BR/-0.60/0.15,
                EL/0.60/0.15,ER/1.70/0.15%
            } {
                \draw[up ancestor] (P)
                    .. controls ({0.55*\tx},2.35)
                        and ({1.05*\tx},{\ty+0.80})
                    .. (\target);
            }

            \foreach \target/\tx in {
                BL/-1.70,BR/-0.60,EL/0.60,ER/1.70%
            } {
                \draw[up ancestor] (A)
                    .. controls ({0.48*\tx},1.60)
                        and ({0.82*\tx},0.90)
                    .. (\target);
            }

            \foreach \u/\v in {
                A/B,A/C,A/D,A/E,B/BL,B/BR,E/EL,E/ER%
            } {
                \draw[up tree] (\u) -- (\v);
            }
            \foreach \target in {BL,BR,EL,ER} {
                \draw[up tree] (\target) -- ++(-0.18,-0.42);
                \draw[up tree] (\target) -- ++(0.18,-0.42);
            }

            \ifnum\stage=2\relax
                \foreach \target/\startangle/\endangle in {
                    B/30/98, C/59/98, D/82/121, E/82/150,
                    BL/59/98, BR/73/98, EL/82/107, ER/82/121%
                } {
                    \draw[up edge mark] (\target)
                        ++(\startangle:0.22)
                        arc[start angle=\startangle,
                            end angle=\endangle,radius=0.22];
                }
            \fi

            \node[up set] at (P) {};
            \foreach \target in {A,B,C,D,E,BL,BR,EL,ER} {
                \node[up vertex] at (\target) {};
            }

            \node[anchor=west] at (0.13,1.84) {$a$};
            \node[anchor=east] at (-1.47,0.74) {$b$};
            \node at (-0.38,0.52) {$c$};
            \node at (0.38,0.52) {$d$};
            \node[anchor=west] at (1.47,0.74) {$e$};

            \foreach \px in {-0.765,0,0.765} {
                \node at (\px,0.80) {$\cdots$};
            }
            \node at (-1.15,0.15) {$\cdots$};
            \node at (1.15,0.15) {$\cdots$};
            \node at (0,-0.98) {($\stage$)};
        \end{scope}
    }
\end{tikzpicture}
\caption{Logical steps in the upward propagation argument at a fixed
second-order stationary point. Each shaded region represents an
internally synchronized set. Starting from the synchronized child
subtrees in (1), the phasor relations illustrated in (2) imply
synchronization of all descendants of $a$ in (3), and then of $a$
itself in (4). The orange arcs in (2) mark only edges from each
indicated node to $a$ and its ancestors, excluding edges within the
corresponding synchronized subtree. The large circle represents $\mathrm{Anc}_T(a)$.
Black edges depict the tree structure, while gray lines represent
additional ancestor--descendant adjacencies in $G$.}
\label{lemma:common-neighbor-sync}
\end{figure}

\subsection{An illustrative example}

We illustrate the upward propagation argument on the graph
$G$ represented by the rooted tree $T$ in
Figure~\ref{proto}, with root $r$.
Fix a second-order stationary point $\boldsymbol{\theta}$
of $E_G$. We first establish synchronization in the subtree
rooted at $c$, then extend the argument to the whole tree.

\begin{figure}[H]
\centering
\begin{tikzpicture}[
    x=1cm, y=1cm,
    font=\normalsize,
    vertex/.style={
        circle, draw=black, fill=white,
        line width=0.65pt,
        minimum size=3.8mm,
        inner sep=0pt, outer sep=0.4pt
    },
    edge/.style={
        line width=0.75pt,
        -{Stealth[length=1.9mm,width=1.35mm]}
    },
    every label/.style={inner sep=1pt},
    label distance=2pt
]
    \node[vertex,label=above:{$r$}] (r) at (0,0) {};
    \node[vertex,label=left:{$b$}]  (b) at (-3.60,-1.35) {};
    \node[vertex,label=right:{$c$}] (c) at (0,-1.35) {};
    \node[vertex,label=right:{$d$}] (d) at (3.60,-1.35) {};

    \node[vertex] (b1)   at (-4.30,-2.70) {};
    \node[vertex] (b11)  at (-4.90,-4.05) {};
    \node[vertex] (b12)  at (-3.70,-4.05) {};
    \node[vertex] (b111) at (-4.90,-5.40) {};

    \node[vertex,label=left:{$c_1$}]  (c1) at (-0.65,-2.70) {};
    \node[vertex,label=right:{$c_5$}] (c5) at (0.65,-2.70) {};
    \node[vertex,label=left:{$c_2$}]  (c2) at (-0.65,-4.05) {};
    \node[vertex,label=below:{$c_3$}] (c3) at (-1.30,-5.40) {};
    \node[vertex,label=below:{$c_4$}] (c4) at (0,-5.40) {};

    \node[vertex] (d1)   at (2.65,-2.70) {};
    \node[vertex] (d2)   at (4.55,-2.70) {};
    \node[vertex] (d11)  at (2.65,-4.05) {};
    \node[vertex] (d21)  at (4.00,-4.05) {};
    \node[vertex] (d22)  at (5.10,-4.05) {};
    \node[vertex] (d111) at (2.00,-5.40) {};
    \node[vertex] (d112) at (3.30,-5.40) {};

    \foreach \parent/\child in {
        r/b, r/c, r/d,
        b/b1, b1/b11, b1/b12, b11/b111,
        c/c1, c/c5, c1/c2, c2/c3, c2/c4,
        d/d1, d/d2, d1/d11, d11/d111, d11/d112,
        d2/d21, d2/d22%
    } {
        \draw[edge] (\child) -- (\parent);
    }
\end{tikzpicture}
\caption{A rooted tree illustrating upward propagation.}
\label{proto}
\end{figure}

We begin with the leaves $c_3$ and $c_4$. Their neighborhoods
in $G$ satisfy
\[
    N_G(c_3)=N_G(c_4)=\{c_2,c_1,c,r\},
\]
so they are structural open twins.
The nonvanishing argument in Step~2 of the proof of
Lemma~\ref{lemma:leaflike}, applied at $c_2$, gives
\[
    \mathbf{v}_{c_2}+\mathbf{v}_{c_1}
    +\mathbf{v}_c+\mathbf{v}_r\neq\boldsymbol{0}.
\]
Thus the phasors of their common neighbors have a nonzero
sum, and Corollary~\ref{stable-open-twins} implies
\[
    \mathbf{v}_{c_3}=\mathbf{v}_{c_4}.
\]

Next, consider the adjacent nodes $c_2$ and $c_3$.
Their common neighbors are $c_1$, $c$, and $r$, while
$c_2$ has the additional neighbor $c_4$:
\[
\begin{aligned}
    N_G(c_2)\setminus\{c_3\}
    &=\{c_1,c,r\}\sqcup\{c_4\},\\
    N_G(c_3)\setminus\{c_2\}
    &=\{c_1,c,r\}.
\end{aligned}
\]
Since $c_4$ synchronizes with $c_3$,
Lemma~\ref{pairwise-no-nauty-extra-node} implies that
$c_2$ also synchronizes with $c_3$. Therefore,
\[
    \mathbf{v}_{c_2}
    =\mathbf{v}_{c_3}
    =\mathbf{v}_{c_4}.
\]

We now move upward to $c_1$. Since $c_2$ is its only child
in $T$, the nodes $c_1$ and $c_2$ have identical closed
neighborhoods in $G$ and are structural closed twins.
By Corollary~\ref{stable-closed-twins}, they synchronize.
Consequently, $c_1$ synchronizes with all its descendants
and is leaf-like in $T$ at $\boldsymbol{\theta}$.

The children of $c$ are $c_1$ and $c_5$.
We have shown that $c_1$ is leaf-like, and $c_5$ is
leaf-like because it is a leaf.
Applying Lemma~\ref{lemma:leaflike} at $c$ therefore shows
that $c$ is leaf-like in $T$ at $\boldsymbol{\theta}$.
Thus all nodes in the subtree rooted at $c$ synchronize.

Repeating this argument in the subtrees rooted at $b$ and
$d$ shows that every child of the root $r$ is leaf-like
in $T$ at $\boldsymbol{\theta}$.
A final application of Lemma~\ref{lemma:leaflike} at $r$
shows that $r$ is also leaf-like in $T$ at
$\boldsymbol{\theta}$.
Since every other node is a descendant of $r$, all nodes
of $G$ synchronize.

\subsection{Proof of Theorem~\ref{main1}: the general case}
\label{general-proof}

\begin{proof}
Let $G=(V,E)$ be a connected quasi-threshold graph with
underlying rooted tree $T$ and root $r$.
Fix an arbitrary second-order stationary point
$\boldsymbol{\theta}$ of $E_G$.

For $v\in V(T)$, let $h(v)$ denote the height of the subtree
rooted at $v$, namely, the maximum number of edges in a path
from $v$ to a leaf in that subtree.
We prove by induction on $h(v)$ that every node $v$ is
leaf-like in $T$ at $\boldsymbol{\theta}$.

If $h(v)=0$, then $v$ is a leaf and
$\mathrm{Desc}_T(v)=\varnothing$.
Thus $v$ is leaf-like in $T$ at $\boldsymbol{\theta}$
by definition.

Suppose that every node of height at most $k$ is leaf-like
in $T$ at $\boldsymbol{\theta}$, and let $u$ be a node
with $h(u)=k+1$.
Then $u$ is a non-leaf node, and every child of $u$
has height at most $k$.
By the induction hypothesis, all children of $u$ are
leaf-like in $T$ at $\boldsymbol{\theta}$.
Lemma~\ref{lemma:leaflike} therefore implies that $u$
is also leaf-like in $T$ at $\boldsymbol{\theta}$.

By induction, every node of $T$ is leaf-like in $T$
at $\boldsymbol{\theta}$.
In particular, the root $r$ is leaf-like, so
\[
    \theta_v\equiv\theta_r\pmod{2\pi}
    \qquad \text{for every } v\in V.
\]
Hence $\boldsymbol{\theta}$ is fully synchronized and
is therefore a global minimizer of $E_G$.
Since $\boldsymbol{\theta}$ was an arbitrary second-order
stationary point, $G$ is second-order globally synchronizing.
\end{proof}